\documentclass[10pt,psamsfonts]{article}
\usepackage{amsmath,amssymb,latexsym}
\usepackage{amscd}
\usepackage[mathscr]{eucal}
\usepackage[all]{xy}
\newcommand\qed{{\hspace*{\fill}$\Box$\vskip12pt plus 1pt}}

\newcommand\se{{\mathscr  E}}
\newcommand\scrf{{\mathscr  F}}

\newcommand\sh{{\mathscr  H}}
\newcommand\sy{{\mathscr I}}
\newcommand\sk{{\mathscr  K}}

\newcommand\sL{{\mathscr  L}}

\newcommand\so{{\mathscr O}}

\newcommand\gra{\alpha}

\newcommand\gre{\varepsilon}

\newcommand\vphi{\varphi}

\newcommand\grb{\beta}
\newcommand\rat{{\mathbb Q}}

\newcommand\zed{{\mathbb Z}}

\newcommand\pn[1]{{\mathbb P}^{#1}}

\newcommand\proof{\noindent{\em Proof.}\ \ }

\def\rank{\mathop{\rm rank}\nolimits}
\def\cd{\mathop{\rm cd}\nolimits}

\def\an{\operatorname{an}}
\def\id{\operatorname{id}}
\def\Spec{\operatorname{Spec}}
\def\Proj{\operatorname{Proj}}

\def\codim{\operatorname{codim}}

\newtheorem{theorem}{Theorem}[section]
\newtheorem{lemma}[theorem]{Lemma}
\newtheorem{corollary}[theorem]{Corollary}

\newtheorem{clam}[theorem]{Claim}
\newtheorem{prop}[theorem]{Proposition}

\newtheorem{deff}[theorem]{Definition}
\newtheorem{re}[theorem]{Remark}

\newtheorem{pargrph}[theorem]{}
\newtheorem{examp}[theorem]{Example}
\newtheorem{examps}[theorem]{Examples}
\def\top{\operatorname{top}}

\newenvironment{rem*}{\begin{re}\em}{\end{re}}
\newenvironment{rems*}{\begin{re}\em}{\end{res}}
\newenvironment{ex*}{\begin{examp}\em}{\end{examp}}
\newenvironment{exs*}{\begin{examps}\em}{\end{examps}}
\newenvironment{deff*}{\begin{deff}\em}{\end{deff}}

\newenvironment{prgrph*}[1]{\indent\begin{pargrph}{\bf #1.}\em\ }{\end{pargrph}}

\begin{document}

\title{Seshadri positive submanifolds of polarized manifolds}

\date{}
\author{Lucian B\u adescu and Mauro C. Beltrametti }
\maketitle
\begin{abstract} Let $Y$ be a submanifold of dimension $y$ of a  polarized complex manifold  $(X,A)$  of dimension $k\geq 3$, with $1\leq y\leq k-1$. We define and study two positivity conditions on $Y$ in $(X,A)$, called Seshadri $A$-bigness and (a stronger one) Seshadri $A$-ampleness. In this way we get the natural generalization of the theory
initiated by Paoletti in \cite{Pao} (which corresponds to the case $(k,y)=(3,1)$) and subsequently generalized and completed in
 \cite{BBF}   (regarding curves in a polarized manifold of arbitrary dimension). The theory presented here, which is new even if $y=k-1$, is motivated by a reasonably large area of examples.
\end{abstract}

\section*{Introduction}

Let  $(X,A)$ be a smooth complex polarized variety (manifold) of dimension $k\geq 3$, and let $Y$ be a smooth (connected) subvariety of dimension $y\geq 1$. Let $X_Y$ be the variety obtained from $X$  by blowing up $Y$, let $\pi\colon X_Y\to X$ be the canonical morphism and let $E=\pi^{-1}(Y)$ be the exceptional divisor of $\pi$. Let $N$ denote the normal bundle of $Y$ in $X$. Then one can define the {\em Seshadri constant} $\varepsilon(Y,A)$ of $Y$ with respect to the polarization $A$ as
$$\varepsilon(Y,A):=\sup\{\eta\in \rat \;|\; A^*-\eta E\;\;{\rm is\;ample}\},$$
where $A^*=\pi^*(A)$. As  $A$ is ample  on $X$ and the line bundle $\mathscr O_{X_Y}(-E)$ is $\pi$-ample (where $\pi\colon X_Y\to X$ is the structural morphism of  $X_Y$), this definition makes sense and yields the inequality $\varepsilon(Y,A)>0$.  Then $\varepsilon(Y,A)$ is a strictly positive real number. Motivated by the study of the gonality of curves and by the behavior  of restriction of the stable vector bundles (see \cite{Pao2}), the Seshadri constant $\varepsilon(Y,A)$ was used by Paoletti in \cite{Pao} to study the so-called Seshadri positive curves in a polarized threefold. Subsequently, this theory was generalized and  completed to the case when $Y$ is a smooth curve in a smooth polarized variety $(X,A)$ of arbitrary dimension $k\geq 3$ in \cite{BBF}. Specifically, for every $\eta\in(0,\varepsilon(Y,A))$ define the  numerical invariant
\begin{equation}\label{def1}\delta_{\eta}(Y,A):=\eta^{k-3}(\eta\deg(N)-(k-2)d),\end{equation}
where $d$ is the degree of $Y$ with respect to the polarization $A$. In this case $Y$ is said to be Seshadri $A$-big in $(X,A)$  if there is  $\eta\in(0,\varepsilon(Y,A))$ such that $\delta_{\eta}(Y,A)>0$. This definition has a natural geometrical interpretation given by Lemma \ref{Big} below.

Moreover, at the end of \cite{BBF} it was left open the problem of finding the natural general setting 
for a theory of Seshadri positivity for a submanifold $Y$ of any dimension $y\geq 1$ in a complex polarized manifold of dimension $k\geq 3$, which should generalize the case of Seshadri positive curves in a polarized manifold. And a possible such definition of Seshadri positivity  was suggested in the last section of \cite{BBF}, but, unfortunately, it turned out to be much too strong to work with if $y\geq 2$.

The aim of the present paper is to give, in our opinion, the natural  generalization of the concept of Seshadri positivity to submanifolds of dimension  $y\geq 1$ in a smooth polarized manifold $(X,A)$ and to recuperate the main results (proved in \cite{Pao} and \cite{BBF} in the case when $Y$ is a curve) in general. To this end, for every submanifold $Y$ of dimension $y\geq 1$ of a polarized manifold $(X,A)$ of dimension $k\geq 3$ and for every $\eta\in(0,\varepsilon(Y,A))$,  define the numerical invariant
\begin{equation}\label{defy}\delta_{\eta}(Y,A):= -\sum_{t=0}^{k-2}{{k-2}\choose{t}}\eta^t\int_Ys_{y-k+t+2}(N)\cdot
c_1(A_Y)^{\cdot (k-t-2)},\end{equation}
where $s_i(N)$ is the $i$-th Segre class  $N$, $c_1(A_Y)$ is the first Chern class of the restriction of $A$ to $Y$, and $\int_Y\alpha$ denotes the degree of a class of a zero-cycle $\alpha$ on $Y$. Thus the coefficients of the polynomial function given by $\eta\mapsto\delta_{\eta}(Y,A)$ depend on the invariants of the closed embedding of $Y$ in the polarized manifold $(X,A)$. It turns out that the polynomial function defined by \eqref{defy}  is the natural generalization of \eqref{def1}. Now, in general, we say that $Y$ is 
Seshadri $A$-big if there exists an $\eta\in(0,\varepsilon(Y,A))$ such that $\delta_{\eta}(Y,A)>0$.
This general definition is geometrically motivated by Lemma \ref{Big} again. Moreover, starting with this definition one can provide a reasonably large area of interesting examples of Seshadri $A$-big manifolds in any dimension and codimension, and one can prove most of the results contained in \cite{Pao} and \cite{BBF} in this general setting. It is worth noting that the theory of Seshadri positivity of a submanifold $Y$ of a polarized manifold $(X,A)$ presented here is of interest even when $Y$ is a hypersurface of $X$ (see Section \ref{codimone}).

The paper is organized as follows. In Section \ref{def} we give the basic definitions of Seshadri $A$-bigness and Seshadri $A$-ampleness for a submanifold $Y$ of dimension 
$y\geq 2$ of a complex polarized projective manifold $(X,A)$ (see Definition \ref{BigAmple}), we prove some preliminary results and we give a few examples of how to compute the Seshadri constant $\varepsilon(Y,A)$. 

In Section \ref{codimone} we discuss the case when $\codim_X(Y)=1$. In this case there are surprisingly many examples. Specifically, we first prove that if the normal bundle $N$ of $Y$ in $X$ is ample then $Y$ is Seshadri $A$-big in $X$. Then we show that the same conclusion holds under a much weaker hypothesis (Corollary \ref{codim111}) and we give three relevant examples in which $N$ is not ample but satisfies that weaker hypothesis.

Section \ref{exampl} is in some sense the ``core" of the paper and deals with the case of submanifolds of codimension $\geq 2$, in which case the concepts of Seshadri $A$-bigness and $A$-ampleness become rather strong. Indeed in this section the task is to  provide interesting examples of submanifolds $Y$ of a polarized manifold $(X,A)$ of codimension $\geq 2$ which are Seshadri $A$-big (in fact, Seshadri $A$-ample by the last statement of Corollary \ref{Same0}). This gives  strong motivations for the positivity notions introduced in Section \ref{def} in arbitrary codimension.  We also point out in paragraph \ref{Sub} a useful expression  which allows us to compute  Chern classes for manifolds in projective space. 
Let us overview the examples we have. As a consequence of a more general result, we prove that,  if we fix a  projective embedding $X\hookrightarrow\mathbb P^n$, every smooth complete intersection in $X$  of codimension $2$ is $\mathscr O_X(1)$-ample (see Propositions \ref{CIcod2} and \ref{Z(s)}). We also show that the  Veronese surface $Y$ in $\pn 5$ is Seshadri $\so_{\pn 5}(1)$-ample, as well as its  projection from a general point of $\pn 5$ is Seshadri $\so_{\pn 4}(1)$-ample. Further relevant examples of subvarieties of $\pn k$ that are Seshadri $\so_{\pn k}(1)$-ample are given by some rational normal scrolls of dimension $2$, and some Segre embeddings (Proposition \ref{F1}). Another example of a Seshadri $A$-big surface in $\mathbb P^5$ is given by the geometrically ruled surface $Y=C \times \pn 1$, with $C\subset\mathbb P^2$ a smooth elliptic curve, embedded in $\pn 5$ via the Segre  embedding $\mathbb P^2\times\mathbb P^1\hookrightarrow\mathbb P^5$ (Proposition \ref{elliptic}).
This latter example has some special interest because the surface in question is not rational. Finally, Proposition \ref{finite}  (noticed in Paoletti \cite{Pao} in the case $(k,y)=(3,1)$) is important because it  allows  to construct many more examples of Seshadri $A$-big (respectively, Seshadri $A$-ample) submanifolds of a polarized manifolds $(X,A)$, starting from some already known ones.

In last two sections we show that the definitions given in Section $2$ and the examples of Section $3$ offer the ``correct setting'' for the study of submanifolds of dimension $\geq 2$ in a polarized manifold $(X,A)$ which are Seshadri positive. In Section \ref{2a} we show how the  Seshadri positivity of  a submanifold  $Y\subset X$ is  related to the theory of formal-rational functions of $X$ along $Y$, as well as to the cohomological dimension of the complement $X\setminus Y$.  Precisely, if  $Y$ is Seshadri $A$-big (respectively, Seshadri $A$-ample), then $Y$ is G$2$ (respectively G$3$) in $X$ in the sense of Hironaka--Matsumura \cite{HM}.  This allows us to give some estimates for the cohomolgical dimension  of  $X\setminus Y$ (see Theorem \ref{Big2}).

In Section \ref{2b} a criterion to distinguish between Seshadri $A$-ampleness and Seshadri $A$-bigness is proved (see Theorem \ref{Main}, compare with \cite[Theorem 3.1]{BBF}).  This criterion asserts that a Seshadri $A$-big subvariety $Y$ is $A$-ample in $X$ if and only if  any irreducible hypersurface in $X$ intersects $Y$. Moreover, we show that if $Y$ is Seshadri $A$-big, the number of all  hypersurfaces of $X$ not meeting $Y$ is finite. As a consequence of Theorem \ref{Main} we can provide  explicit examples of  submanifolds $Y\subset X$ which are Seshadri $A$-big, but not Seshadri $A$-ample (see Example \ref{exceptional}).
The section ends up with some further examples, comments and remarks.

Unless otherwise stated, the terminology and the notation used throughout are standard. As far as the Chern and the Segre classes computations are concerned (see Sections \ref{def} and \ref{exampl}) we shall follow the notation, conventions and some basic results of Fulton's book \cite{Fu}. All algebraic varieties (or schemes) will be defined over the field $\mathbb C$ of complex numbers.

\smallskip

{\small{\bf Dedication.}  In the middle nineties (inspired by \cite{Pao}) we started studying these kind of problems  for curves in a polarized manifold of arbitrary dimension in \cite{BBF}, together  with our colleague and friend Paolo Francia (who passed away in July 2000). We want to dedicate this paper to his memory.}

\section{Basic  definitions and preliminary results and examples}
\label{def}
\addtocounter{subsection}{1}\setcounter{theorem}{0}

Let $X$ be a smooth complex manifold of dimension $k\geq 3$ and let $Y\subset
X$ be a smooth  subvariety of $X$ of dimension $y$ such that $1\leq y\leq
k-2$. Let $N:=N_{Y|X}$ be the normal bundle of $Y$ in $X$. Then $N$ is a
rank $k-y$ vector bundle on $Y$.  Let
\begin{equation}\label{diag}
\xymatrix{E \  \ar@{^(->}[r]^{j}  \ar[d]_{\pi'} & X_Y \ar[d]^{\pi} \\Y \ \ar@{^(->}[r]^{i} & X }
 \end{equation}
be the blowing up of $X$ along $Y$,  with $E:=\pi^{-1}(Y)$ the exceptional divisor and $i$ and $j$ the natural closed immersions.
Then $E={\mathbb P}(N^{\vee})$, where $N^{\vee}$ is the dual of $N$. We are using the notation of Grothendieck, i.e., if $\mathscr F$ is a vector bundle on an algebraic  variety $Z$, we denote by 
$\mathbb P(\mathscr F)$    the projective bundle of hyperplanes of $\mathscr F$ and by $\mathscr O_{\mathbb P(\mathscr F)}(1)$ is the tautological line bundle of $\mathbb P(\mathscr F)$.
For any ample line bundle $A$ on $X$, set $A^*:=\pi^*(A)$ and let
$d:=Y\cdot A^{\cdot y}$ be the {\em degree} of $Y$ with respect to $A$.

According to Fulton \cite[(3.1), (B.5.5)]{Fu}, if $\se$ is  a vector bundle on an algebraic variety $Z$, define the {\em Segre class} $s_j(\se)$ of order $j$ of $\se$ by
\begin{equation}\label{segre0}s_j(\se):=\pi_*(c_1(\so_{{\mathbb P}(\se^{\vee})}(1))^{\cdot{\rm (rank}(\se)-1+j)}),\end{equation}
where $\pi\colon\mathbb P(\se^{\vee})\to Z$ is the canonical projection. Then $s_j(\se)$ is the rational equivalence class of a cycle of codimension $j$ in $Z$. In particular, $s_0(\se)=1$ and $s_j(\se)=0$ if $j<0$ or  $j>\dim(Z)$. Moreover, if $\rank(\se)=1$ then 
$\mathbb P(\se^{\vee})=\mathbb P(\se^{-1})=Z$, $\pi=\id_Z$ and $\mathscr O_{\mathbb P(\se^{\vee})}(1)=\se^{-1}$. In particular, if $\se$ is invertible it follows that $s_j(\se)=(-1)^jc_1(\se)^{\cdot j}$ for every $j\geq 0$.

First we need the following:

\begin{lemma}\label{relations} With the  notation as in diagram $(\ref{diag})$, set $A_Y:=A\otimes\mathscr O_Y$. Then for every nonnegative integer  $r$  we have the following equalities:
$$ \int_{X_Y} E^{\cdot (k-r)}\cdot c_1(A^*)^{\cdot r}=\begin{cases} (-1)^{k-r-1} \int_Ys_{y-r}(N)\cdot c_1(A_Y)^{\cdot r} & \text{if} \ r\leq y,\\ 
0 &\text{if} \  r\geq y+1.\end{cases}$$
\end{lemma}

\proof  One has $\so_{X_Y}(-E)_E=\so_{{\mathbb P}(N^{\vee})}(1)=\so_E(1)$.
Therefore we get the following equalities of $r$-cycles (modulo the rational equivalence)
\begin{eqnarray*}
E^{\cdot (k-r)} &=&  E_E\cdots E_E=c_1(\so_E(-1))\cdots c_1(\so_E(-1))\;\;\;{\small (k-r-1 \;{\rm
times}})\\&=&  (-1)^{k-r-1}c_1(\so_{{\mathbb P}(N^{\vee})}(1))^{\cdot (k-r-1)}.
\end{eqnarray*}
Thus,  by using the projection formula (\cite[Proposition 2.5, (c)]{Fu}) and expression \eqref{segre0}, we get
\begin{eqnarray*}
\int_{X_Y}E^{\cdot (k-r)}\cdot c_1(A^*)^{\cdot r} & = &(-1)^{k-r-1}\int_{X_Y}c_1(\so_{{\mathbb P}(N^{\vee})}(1))^{\cdot (k-r-1)}\cdot
c_1(A^*)^{\cdot r} \\  &=&(-1)^{k-r-1}\int_X\pi_*(c_1(\so_{{\mathbb P}(N^{\vee})}(1))^{\cdot (k-r-1)})\cdot c_1(A)^{\cdot r}\\&=&(-1)^{k-r-1}\int_Y\pi'_*(\so_{{\mathbb P}(N^{\vee})}(1))^{\cdot (k-r-1)})\cdot c_1(A_Y)^{\cdot r}\\&=&(-1)^{k-r-1}\int_Ys_{y-r}(N)\cdot c_1(A_Y)^{\cdot r}.\end{eqnarray*}
Finally, if $r\geq y+1$ then $y-r<0$, whence $s_{y-r}(N)=0$. \qed

\begin{deff*}(Paoletti \cite{Pao})\label{Seshadri}  Let $X$ be a smooth complex manifold of dimension $k\geq 3$,  let $Y\subset X$ be a smooth  subvariety of $X$ of positive dimension $\dim (Y):=y\leq k-2$, and fix an ample line bundle $A$ on $X$. Define the {\em Seshadri constant} of $\,Y$ with respect to  $A$  as
\begin{equation}\label{Ses}
\gre(Y,A):=\sup\{\eta\in \rat \;|\; A^*-\eta E\;\;{\rm is\;ample}\}.
\end{equation}

We claim that $0<\varepsilon(Y,A)<\infty$.  Indeed, since $A$ is ample on $X$ and $\mathscr O_{X_Y}(-E)$ is $\pi$-ample on $X_Y$, by \cite[II (4.6.13), (ii)]{EGA}, we know that $A^*-\frac{1}{n}E$ is ample for $n\gg 0$, whence $\varepsilon(Y,A)>0$. On the other hand,
for every ample divisor $H$ on $X_Y$ one has $H^{k-1}\cdot
(A^*-\eta E)<0$ for $\eta \gg 0$, so that $\varepsilon(Y,A)<\infty$.

Define also the numerical invariant $\delta(Y,A)$  by
\begin{equation}\label{delta2}
\delta(Y,A):=
-\sum_{t=0}^{k-2}{{k-2}\choose{t}}\varepsilon(Y,A)^t\int_Ys_{y-k+t+2}(N)\cdot
c_1(A_Y)^{\cdot (k-t-2)}.\end{equation}
Moreover, for every $\eta\in(0,\varepsilon(Y,A))$ also define 
\begin{equation}\label{delta1}
\delta_{\eta}(Y,A):= -\sum_{t=0}^{k-2}{{k-2}\choose{t}}\eta^t\int_Ys_{y-k+t+2}(N)\cdot
c_1(A_Y)^{\cdot (k-t-2)}.\end{equation}
\end{deff*}

The numerical invariants  $\varepsilon(Y,A)$, $\delta(Y,A)$ and $\delta_{\eta}(Y,A)$ will play an  important role in the studying  the embedding $i\colon Y\hookrightarrow X$. For instance the numerical invariant $\delta_{\eta}(Y,A)$ will be motivated by Lemma \ref{Big} below.
From the definitions it follows immediately  that  
\begin{equation}\label{equivalence}
\delta(Y,A)>0 \Longrightarrow\delta_{\eta}(Y,A)>0\;\text{for every $\eta\in(0,\varepsilon(Y,A))$ which is close to $\varepsilon(Y,A)$.}
\end{equation}

Moreover, if $\delta(Y,A)=0$ and if $f'(\varepsilon(Y,A))<0$ then $\delta_{\eta}(Y,A)>0$ for every $\eta\in(0,\varepsilon(Y,A))$ which is close to $\varepsilon(Y,A)$, where $f\colon\mathbb R\to\mathbb R$ is the polynomial function defined by $f(\eta)=\delta_{\eta}(Y,A)$ for all $\eta\in\mathbb R$, and $f'$ denotes the derivative of $f$.

\begin{rem*}\label{indep} From the definition of the Seshadri constant it follows that $\varepsilon(Y,A^{\otimes s})=s\varepsilon(Y,A)$, for every integer $s\geq 1$. In particular, for every $\eta\in(0,\varepsilon(Y,A))$ we have $s\eta\in(0,\varepsilon(Y,A^{\otimes s}))$. Moreover, from  definition \eqref{delta1} of $\delta_{\eta}(Y,A)$ it follows that 
\begin{equation}\label{indep1}\delta_{s\eta}(Y,A^{\otimes s})=s^{k-2}\delta_{\eta}(Y,A),\;\;\text{for every $\eta\in(0,\varepsilon(Y,A))$ and $s\geq 1$.}\end{equation}
\end{rem*}

\begin{ex*}\label{Remeq} 
In the special case $y=1$ (studied in \cite{BBF}) the converse of (\ref{equivalence}) is also true, i.e., 
$$\delta_{\eta}(Y,A)>0\;\;\text{for some $\eta\in(0,\varepsilon(Y,A))$}\;\Longrightarrow\;\delta(Y,A)>0.$$
Indeed in this case  $y-k+t+2= -k+t+3$, whence $s_{-k+t+3}(N)=0$ for $t<k-3$. As $s_1(N)=-\deg N$,
\begin{equation}\label{delta4} \delta(Y,A)=\varepsilon(Y,A)^{k-3}(\varepsilon(Y,A)\deg N-(k-2)d)\;\,\text{and}\;\,
\delta_{\eta}(Y,A)=\eta^{k-3}(\eta\deg N-(k-2)d), \end{equation}
where $d:=\int_Yc_1(A_Y)$ is the degree of $Y$ with respect to the polarization $A$.
In particular we recover  the definitions of $\delta(Y,A)$ and $\delta_{\eta}(Y,A)$ given in \cite{Pao} (for $k=3$) and in \cite{BBF} (for $k\geq 4$) when $Y$ is a curve. Moreover, in this case, we also get that
$\delta_\eta(Y,A)>0$ implies $\deg N>0$, whence
$$\delta(Y,A)=\gre^{k-3}(\gre\deg N-(k-2)d)>\eta^{k-3}(\eta\deg N-(k-2)d)=\delta_{\eta}(Y,A)>0,$$
since $\gre:=\gre(Y,A)>\eta$.
\end{ex*}

\begin{ex*}\label{Remeq2} If $Y$ is a surface  formulae \eqref{delta2}, \eqref{delta1} become
\begin{equation*} \begin{split}
 \delta(Y,A)&=-\varepsilon^{k-4}\Big(\varepsilon^2\! \int_Ys_2(N)-\varepsilon (k-2)\!\! \int_Ys_1(N)\cdot c_1(A_Y)-\frac{(k-2)(k-3)}{2}\int_Yc_1(A_Y)^{\cdot 2}\Big), \\
 \delta_{\eta}(Y,A)&=-\eta^{k-4}\Big(\eta^2\!\int_Ys_2(N)-\eta(k-2)\!\! \int_Ys_1(N)\cdot c_1(A_Y)-\frac{(k-2)(k-3)}{2}\int_Yc_1(A_Y)^{\cdot 2}\Big),
\end{split} \end{equation*} 
where $\varepsilon:=\varepsilon(Y,A)$ and $\eta\in(0,\varepsilon)$.
From the equalities  $s_1(N)=-c_1(N)$ and $s_2(N)=c_1(N)^{\cdot 2}-c_2(N)$ (which follow from the definition of Segre and Chern classes, see Example \ref{codim2} below) we get
\begin{equation} 
\delta = -\varepsilon^{k-4}\Big(\varepsilon^{2}\!\! \! \int_Y(c_1(N)^{\cdot 2}-c_2(N))-\varepsilon (k-2)\!\!\!\int_Yc_1(N)\cdot c_1(A_Y)+\frac{(k-2)(k-3)}{2}\deg(Y)\Big), \label{Ses2} \end{equation}  and
\begin{equation} 
 \delta_{\eta} =  -\eta^{k-4}\Big(\eta^2\!\! \!\int_Y(c_1(N)^{\cdot 2}-c_2(N))-\eta(k-2) \!\!\! \int_Yc_1(N)\cdot c_1(A_Y)+\frac{(k-2)(k-3)}{2}\deg(Y)\Big),\label{Ses2'}
\end{equation} 
where $\delta:=\delta(Y,A)$,  $\delta_{\eta} := \delta_{\eta}(Y,A)$, and  $\deg(Y)=\int_Yc_1(A_Y)^{\cdot 2}$ is the degree of $Y$ with respect to the polarization $A$.
\end{ex*}

Here are a couple of examples regarding the computation of Seshadri constants.

\begin{ex*}\label{CI}(Complete intersections)  Let $X$ be a submanifold of dimension $k$ of the projective space $\mathbb P^n$, and let $Y$ be a 
submanifold of dimension $y$ of $X$. Assume that $Y$  is a (scheme
theoretic) complete intersection in $X$, i.e., $Y$  is the transversal
intersection $Y=F_1\cap\cdots\cap F_{k-y}$  of $k-y$  hypersurfaces
$F_i\subset X$, $i=1,\ldots,k-y$, of $\mathbb P^n$.  Let $A=\so_X(1)$.
Let $d_i={\rm deg}F_i$, $i=1,\ldots,k-y$, and assume $d_1\geq\cdots\geq
d_{k-y}$. Then
the Seshadri constant of $(Y,A)$ is $$\gre(Y,A)=\frac{1}{d_1}.$$
It is a general well known fact that $d_1A^*-E$ is spanned by its global
sections, which implies that 
$\gre(Y,A)\geq\frac{1}{d_1}$ (compare with  Example \ref{MinDeg} below).

Therefore, to prove the equality, it is enough to show that $d_1A^*-E$ is not
ample. To see this note that the restriction of $d_1A^*-E$ to the section
of $E\to Y$
 corresponding to the quotient $N^{\vee}\to (-d_1A^*)\to 0$ is the trivial bundle (see also \cite[Corollary 2.5]{BS}).
\end{ex*}

\begin{ex*}\label{MinDeg} (Varieties defined in a given degree) Notation as in Definition \ref{Seshadri}. Let $\sy_Y$ be the ideal sheaf of $Y$ in $X$. We say that $Y$ is {\em defined by $A$ in degree} $t$ if either $\sy_Y(t):=A^{\otimes t}\otimes \sy_Y$ is generated by its global sections, or, equivalently, if $Y$ is the scheme theoretic intersection of all divisors ${\rm div}_X(s)$, with $s\in H^0(X, A^{\otimes t}\otimes \sy_Y)\setminus\{0\}$.
In this case, from  \cite[Lemma 1.5 and Corollary 1.6]{BBF}, we know that 
$$\gre (Y,A)\geq\frac{1}{t}.$$

The degree $t=2$ case with $X=\pn k$, $Y$ a submanifold of $\pn k$ (not a linear subspace of $\pn k$) and $A=\so_{\pn k}(1)$ is of special interest in the sequel. Take a line $\ell$ meeting $Y$ in at least two points but not contained in $Y$. Note that such a line exists, since otherwise $Y$ would be a linear subspace in $\pn k$. Let $\ell'$ be the proper transform of  $\ell$ under  the blowing up $\pi$ of $\pn k$ along $Y$. Then $$\pi^*\sy_Y(2)\cdot \ell'=(2A^*-E)\cdot \ell'=2-E\cdot\ell'\leq 0.$$ Therefore 
$2A^*-E$ is spanned but not ample. It thus follows  that 
\begin{equation}\label{12} \gre (Y,A)=\frac{1}{2}.\end{equation}

Many of the varieties $Y\subset \pn k$ we will deal with throughout the paper---varieties of minimal degree, Segre varieties---are indeed generated in degree two, so that they satisfy condition (\ref{12}). Determinantal varieties give further examples of varieties defined in a given degree (see \cite[\S 2]{BS}).

\end{ex*}

\begin{rem*}\label{DDD} As $(0,\varepsilon(Y,A)\cap\mathbb Q)$ is dense in $(0,\varepsilon(Y,A))$, it follows that there exists $\eta\in(0,\varepsilon(Y,A))$ such that $\delta_{\eta}(Y,A)>0$ if and only if there exists $\eta\in(0,\varepsilon(Y,A))\cap\mathbb Q$ such that $\delta_{\eta}(Y,A)>0$.\end{rem*}

For the simplicity of notation, when there is no danger of confusion, we shall make no difference between the line bundle $A^*$ and the Chern class $c_1(A^*)$. Choose now $\eta\in (0,\gre(Y,A))\cap\mathbb Q$ such that the $\mathbb Q$-line bundle $A^*\otimes\mathscr O_{X_Y}(-\eta E)=\mathscr O_{X_Y}(A^*-\eta E)$ is ample. Write
$\eta=\frac{n}{m}$, with $m$ and $n$ positive integers.  Then the divisor $mA^*-nE$ is ample by the choice of $\eta=\frac{n}{m}$. Assume that in the ample linear system $|mA^*-nE|$ there exist divisors $H_1,\ldots,H_{k-2}$ such that  $\mathscr H$ is a smooth surface and $Y'$ is a smooth curve, where
\begin{equation}\label{Y'}Y':=E\cap \sh, \;\;{\rm and}\; \;\sh:=H_1\cap \cdots\cap H_{k-2}.\end{equation}
 Note that this condition can be easily realized because if we
multiply $m$, $n$ by a suitable positive large integer $\lambda$ (which does not change $\eta$), the linear system $|\lambda mA^*-\lambda nE|$ becomes very ample, and then apply Bertini's theorem.

\begin{lemma}\label{Big} Let $\eta=\frac{n}{m}\in (0,\gre(Y,A))\cap\mathbb Q$, with $n$ and $m$  positive integers. Then $\delta_{\eta}(Y,A)>0$ if and
only if $(Y'^{\cdot 2})_\sh>0$, where $Y'$ and $\mathscr H$ are defined by \eqref{Y'} and $(Y'^{\cdot 2})_{\mathscr H}:=\int_{\mathscr H}Y'^{\cdot 2}$.\end{lemma}

\proof As
$(Y'^{\cdot2})_{\sh}=\int_{\mathscr H}E_{\sh}^{\cdot 2}=\int_{X_Y}E^{\cdot2}\cdot(mc_1(A^*)-nE)^{\cdot(k-2)}$, by Lemma \ref{relations} we get
\begin{eqnarray*}\frac{1}{m^{k-2}}(Y'^{\cdot2})_{\sh} &=&\int_{X_Y}E^{\cdot2}\cdot(c_1(A^*)-\eta E)^{\cdot(k-2)} \\
&=& \sum_{t=0}^{k-2}(-1)^t\eta^t{{k-2}\choose{t}}\int_{X_Y} E^{\cdot(t+2)}\cdot {c_1(A^*)}^{\cdot(k-2-t)}\end{eqnarray*}
\begin{eqnarray*}&=&
\sum_{t=0}^{k-2}(-1)^t{{k-2}\choose{t}}\eta^t(-1)^{t+1}\int_Ys_{y-k+t+2}(N)\cdot
c_1(A_Y)^{\cdot(k-t-2)} \\&=&
-\sum_{t=0}^{k-2}{{k-2}\choose{t}}\eta^t\int_Ys_{y-k+t+2}(N)\cdot
c_1(A_Y)^{\cdot(k-t-2)}\\&=&\delta_\eta(Y,A).\end{eqnarray*}\vskip -0.6cm\qed

\vskip 0.2cm

Now we come up to the following basic general definition.

\begin{deff*}\label{BigAmple} Under the above notation, we say that the submanifold $Y$ of $X$ is {\em Seshadri $A$-big} (or {\em big with respect to a fixed polarization $A$ of $X$}) if there exists an $\eta\in(0,\varepsilon(Y,A))$ such that $\delta_{\eta}(Y,A)>0$. By Remark \ref{DDD} the latter condition is also equivalent to the existence of an $\eta'\in(0,\varepsilon(Y,A))\cap\mathbb Q$ such that $\delta_{\eta'}(Y,A)>0$. Or more geometrically, in view of Lemma \ref{BigAmple}, the submanifold $Y$ of $X$ is 
Seshadri $A$-big  if and only if there exists $\eta=\frac{n}{m}\in (0,\gre(Y,A))\cap\mathbb Q$ such that $(Y'^{\cdot 2})_{\sh}>0$. As $Y'$ is a smooth irreducible curve on the smooth projective surface $\mathscr H$, the fact that 
$Y\subset X$ is Seshadri $A$-big amounts to the fact that the normal bundle $N_{Y'|\mathscr H}=\so_\sh(Y')\otimes\so_{Y'}$ of $Y'$ in $\mathscr H$ is ample on $Y'$. 

We also say
that $Y\subset X$ is {\em Seshadri $A$-ample} (or {\em ample with respect to the polarization $A$}) if $Y'$ is an ample divisor on $\sh$. Clearly, if $Y$ is Seshadri $A$-ample then $Y$ is Seshadri $A$-big. The  Seshadri $A$-ampleness is also an open condition, i.e., if it is satisfied for some
$\eta=\frac{n}{m}\in (0,\gre(Y,A))\cap\mathbb Q$, then it also satisfied for every $\eta'\in(\eta-\nu,\eta+\nu)\cap \mathbb Q$, with $\nu>0$ sufficiently small.\end{deff*}

\begin{rem*} From Definition \ref{BigAmple} and  \eqref{indep1} it follows immediately that a submanifold $Y$ of $X$ is Seshadri $A$-big if and only if it is Seshadri $A^{\otimes s}$-big for each integer $s\geq 1$. Using this and Theorem \ref{Main} below we shall see that a submanifold $Y$ of $X$ is Seshadri $A$-ample if and only if it is Seshadri $A^{\otimes s}$-ample.
\end{rem*}
\begin{prgrph*}{Some generalities on Segre and Chern classes}\label{codim2} It is usually more convenient  to express the numerical invariants $\delta(Y,A)$ and $\delta_{\eta}(Y,A)$ in terms of the Chern classes of the normal  bundle $N$ 
of a  closed smooth $y$-dimensional subvariety $Y$ of the projective smooth  $k$-dimensional variety $X$. 
This is due to the fact that in general a vector bundle $\se$ of rank $r$ on an algebraic variety $Z$ has Chern classes  $c_i(\se)=0$ for every $i>\min\{\dim(Z),r\}$. For this we recall the identity
$s_t(\se)\cdot c_t(\se)=1$ between the Segre and the Chern polynomials $s_t(\se)$ and $c_t(\se)$, where
$$s_t(\se)=\sum_{i=0}^{y}s_i(\se)t^i\;\;\text{and}\;\;c_t(\se)=\sum_{i=0}^{y}c_i(\se)t^i,\;\;\text{with $s_0(\se)=c_0(\se)=1$,}$$
which  holds for every vector bundle $\se$. (See  \cite[p. 50]{Fu} for the definition of the Chern polynomials.)
And, for every $n\geq 1$,  we have the  general recurrence formula 
$$s_n(\se)=-s_{n-1}(\se)\cdot c_1(\se)-s_{n-2}(\se)\cdot c_2(\se)-\cdots-s_1(\se)\cdot c_{n-1}(\se)-c_n(\se).$$
In particular,
\begin{equation}\label{sc1}
s_1(\se)=-c_1(\se), \;\; s_2(\se)=c_1(\se)^{\cdot 2}-c_2(\se), \;\;
s_3(\se)=2c_1(\se)\cdot c_2(\se)-c_1(\se)^{\cdot 3}-c_3(\se), \end{equation}
\begin{equation}\label{sc2}s_4(\se)=-3c_1(\se)^{\cdot 2}\cdot c_2(\se)+c_1(\se)^{\cdot4}+c_2(\se)^{\cdot 2}+2c_1(\se)\cdot c_3(\se)-c_4(\se),\ldots \end{equation}
 Recall also that $c_1(\se)=c_1(\det(\se))$, where $\det(\se)=\wedge^r\se$, see \cite[p. 55]{Fu}. 

Take now $Z=Y$ and $\se=N$ (the normal bundle of $Y$ in $X$). If for example
$\codim_X(Y)=2$ and $k\leq 6$ we have $c_i(N)=0$ for every $i\geq 3$, whence  the above formulae and \eqref{delta1} yield  for $\delta_{\eta}(Y,A)$ the following expressions:
\begin{eqnarray}
 \delta_{\eta}(Y,A)&=&-\eta^2\int_Y(c_1(N)^{\cdot 2}-c_2(N))\nonumber \\ & &+2\eta\int_Yc_1(N)\cdot c_1(A_Y)   -\deg(Y), \;\;\text{for $(k,y)=(4,2)$};\label{y=2}
 \end{eqnarray}
\begin{eqnarray}
 \delta_{\eta}(Y,A)&=&-\eta^3\int_Y(2c_1(N)\cdot c_2(N))-c_1(N)^{\cdot 3}) \nonumber \\ & & -3\eta^2\int_Y(c_1(N)^{\cdot 2}-c_2(N))\cdot c_1(A_Y) \nonumber \\
 & & +3\eta\int_Yc_1(N)\cdot c_1(A_Y)^{\cdot2}-\deg(Y), \;\;\text{for $(k,y)=(5,3)$};\label{y=3}
\end{eqnarray}
\begin{eqnarray}
  \delta_{\eta}(Y,A)&=&\eta^4\int_Y(3c_1(N)^{\cdot2}\cdot c_2(N)-c_2(N)^{\cdot2}-c_1(N)^{\cdot4}) \nonumber \\
& & -4\eta^3\int_Y(2c_1(N)\cdot c_2(N)-
c_1(N)^{\cdot3})\cdot c_1(A_Y) \nonumber \\
 & & -6\eta^2\int_Y(c_1(N)^{\cdot 2}-c_2(N))\cdot c_1(A_Y)^{\cdot2}
\nonumber \\ & & +4\eta\int_Yc_1(N)\cdot c_1(A_Y)^{\cdot3}-\deg(Y), \;\;\text{for $(k,y)=(6,4)$}.\label{y=4}
\end{eqnarray}\end{prgrph*}

\section{The one-codimensional case}\label{codimone}
\addtocounter{subsection}{1}\setcounter{theorem}{0}

 Let us start with the simplest case when $\codim_X(Y)=1$, i.e., $y=k-1$. In this case, under the notation of diagram \eqref{diag}, we have $X_Y=X$ and $E=Y$. Therefore 

\begin{equation}\label{Y1}\varepsilon(Y,A)=\sup\{\eta\in \rat \;|\; A-\eta Y\;\;{\rm is\;ample}\}.\end{equation}

The following result already provides    many examples of Seshadri $A$-big submanifolds of codimension one.

\begin{prop}\label{codim1} Let $(X,A)$ be a polarized manifold of dimension $k\geq 2$, and let $Y$ be a submanifold of $X$ of dimension $y=k-1$. If the normal bundle $N=\mathscr O_X(Y)\otimes \mathscr O_Y$ of $Y$ in $X$ is ample then $Y$ is Seshadri $A$-big. If $k=2$ the converse also holds. \end{prop}
\proof  Since $y=k-1$,  formula \eqref{delta1} yields
\begin{eqnarray*}
\delta_{\eta}(Y,A) & = &-\sum_{t=0}^{k-2}\binom{k-2}{t}\eta^t\int_Ys_{t+1}(N)\cdot c_1(A_Y)^{\cdot(k-t-2)} \\  &=&-\sum_{t=0}^{k-2}\binom{k-2}{t}\eta^t(-1)^{t+1}\int_Yc_1(N)^{\cdot (t+1)}\cdot c_1(A_Y)^{\cdot(k-t-2)}\\&=& \int_Yc_1(N)\cdot c_1(A_Y-\eta N)^{\cdot(k-2)}.\end{eqnarray*}
Thus $\delta_\eta(Y,A)>0$
because $N$ and $A_Y-\eta N$ are ample on $Y$  for every $\eta\in(0,\varepsilon(Y,A))$ (the latter as the restriction of the ample line bundle $A-\eta Y$, via \eqref{Y1}). Therefore $Y$ is Seshadri $A$-big.

Conversely, assume $k=2$ and $\delta_{\eta}(Y,A)>0$. Then $\deg(N)=\int_Yc_1(N)=\delta_{\eta}(Y,A)>0$, whence $N$ is ample because it is a line bundle on the curve $Y$.
\qed

In order to provide more examples of Seshadri $A$-big submanifolds in codimension one we need the following  straightforward consequence of Lemma \ref{Big} and Definition \ref{BigAmple}.

\begin{lemma}\label{codim11} Let $(X,A)$ be a polarized manifold of dimension $k\geq 2$, and let $Y$ be a submanifold of $X$ of dimension $y=k-1$. Then 
the following conditions are equivalent:
\begin{enumerate}
\em\item\em $Y$ is Seshadri $A$-big. 
\em\item\em There exists an $\eta=\frac{n}{m}\in(0,\varepsilon(Y,A))\cap\mathbb Q$, with $m$ and $n$ positive integers, such that the line bundle $mA-nY$ is very ample on $X$ and for every general divisors  $H_1,\ldots, H_{k-2}\in|mA-nY|$ one has $(Y'^{\cdot 2})_{\mathscr H}:=\int_{\mathscr H}Y'^{\cdot 2}>0$, where $Y':=Y\cap H_1\cap\cdots\cap H_{k-2}$ and
$\mathscr H:=H_1\cap\cdots\cap H_{k-2}$ $($see \eqref{Y'}$)$.
\em\item\em There exist an $\eta=\frac{n}{m}\in(0,\varepsilon(Y,A))\cap\mathbb Q$, with $m$ and $n$ positive integers,
and divisors $H_1,\ldots, H_{k-2}\in|mA-nY|$ such that $\mathscr H:=H_1\cap\cdots\cap H_{k-2}$ is a smooth surface,  $Y':=Y\cap H_1\cap\cdots\cap H_{k-2}$ a smooth curve, and  $(Y'^2)_{\mathscr H}:=\int_{\mathscr H}Y'^{\cdot 2}>0$.
\end{enumerate}
\end{lemma}

\begin{corollary}\label{codim111}  Let $(X,A)$ be a polarized manifold of dimension $k\geq 3$, and let $Y$ be a submanifold of $X$ of dimension $y=k-1$. Assume that  the normal bundle $N=\mathscr O_X(Y)\otimes \mathscr O_Y$ of $Y$ in $X$ satisfies the following condition:
\begin{equation}\label{22} \deg(N\otimes\mathscr O_C))>0\;\;\text{for every smooth curve $C$ in $Y$ with normal bundle $N_{C|Y}$ ample.}\end{equation}
 Then $Y$ is Seshadri $A$-big.
\end{corollary}

\proof  It is enough to check  condition $2)$ of Lemma  \ref{codim11}. Recall that  for $\eta=\frac{n}{m}\in(0,\varepsilon(Y,A))\cap\mathbb Q$, with $m$ and $n$ positive integers such that the line bundle $mA-nY$ is very ample on $X$, and for every general divisors  $H_1,\ldots, H_{k-2}\in|mA-nY|$, the intersection $Y'=Y\cap\mathscr H$ is proper, where $\mathscr H=H_1\cap \cdots\cap H_{k-2}$ is a smooth surface. Therefore  $N_{Y'|\mathscr H}\cong N\otimes\mathscr O_{Y'}$, and hence $N_{Y'|\mathscr H}$ is ample. It follows that  $(Y'^{\cdot 2})_{\mathscr H}=\deg(N_{Y'|\mathscr H})=\deg(N\otimes\mathscr O_{Y'})>0$ by hypothesis \eqref{22}.
\qed

\begin{exs*}\label{Mumford} Corollary \ref{codim111} leads to explicit construction of $A$-big divisors.

i) In characteristic zero there are examples of line bundles $N$ on some smooth projective varieties $Y$ of dimension $\geq 2$ such that 
$\deg(N\otimes\mathscr O_C)>0$ for every irreducible curve $C$ on $Y$, but which are not ample. The first of such
 examples has been found by Mumford (see e.g. \cite[Example 10.6, p. 56]{Ample}). Specifically, let $B$ be a smooth projective curve of genus $g\geq 2$ over ${\mathbb C}$  and let $E$ be a vector bundle of rank $2$  and degree $0$ over $B$ such that all symmetric powers 
${\bf S}^m(E)$, $m\geq 1$, are stable (one can show that such vector bundles $E$ do exist). Let $Y=\mathbb P(E)$ be the projective bundle associated  to $E$ and set $N=\mathscr O_{\mathbb P(E)}(1)$. Then $\int_Yc_1(N)^{\cdot 2}=\deg(E)=0$ (in particular  $N$ is not ample), but one can show that condition \eqref{22} 
of Corollary \ref{codim111} is satisfied  \cite[Example 10.6, p. 56]{Ample}. Moreover in this case one has $H^0(Y,N)=0$.

Starting with Mumford's example $Y=\mathbb P(E)$ above, we may easily construct an example of a polarized threefold $(X,A)$ containing $Y$ as a divisor with normal bundle $N_{Y|X}=N=\mathscr O_{\mathbb  P(E)}(1)$ not ample, but satisfying condition \eqref{22} of Corollary 
\ref{codim111}. 
Indeed, let  $X':=\mathbb V(N^{-1})=\underline{\Spec}({\bf S}(N^{-1}))$ be the vector bundle over $Y$ associated to the dual of $N$ (with ${\bf S}(N^{-1})=\oplus_{m=0}^{\infty}N^{\otimes{-m}}$) and let 
$i\colon Y\hookrightarrow X'$ be the zero section of $\mathbb V(N^{-1})$ (which corresponds to the canonical surjection $\oplus_{m=0}^{\infty}N^{\otimes{-m}}\to\mathscr O_Y$). Then the normal bundle of $i(Y)$ in $X'$ is isomorphic to $N$. Since $X'$ is smooth and only quasi-projective, take the natural projective closure $X:=\mathbb P(N^{-1}\oplus\mathscr O_Y)$ of $X':=\mathbb V(N^{-1})$). Finally, choosing any polarization $A$
on $X$ we get an example of a polarized threefold $(X,A)$  containing $i(Y)\cong Y$ with normal bundle (isomorphic to) $N$ which is not ample, but it satisfies the hypotheses of Corollary \ref{codim111}. Then  $i(Y)$ is Seshadri $A$ big in $X$.

\smallskip

ii) Note also that there is an example (due to C.P. Ramanujam \cite[Example 10.8, p. 57]{Ample}) of a projective threefold $Y$ over $\mathbb C$ and a line bundle $N$ on $Y$ which is not ample, but still satisfies condition \eqref{22} of Corollary \ref{codim111}. Moreover, in Ramanujam's example one has $\int_Yc_1(N)^{\cdot 3}>0$ and  $H^0(Y,N)\neq 0$. Then, exactly as in case ii)  above, starting with the threefold $Y$ and with the line bundle $N$, one can construct an example of a polarized fourfold $(X,A)$ containing $Y$ as a divisor with normal bundle $N$ and satisfying the hypotheses of Corollary \ref{codim111} (taking $X:=\mathbb P(N^{-1}\oplus\mathscr O_Y)$ and for $A$ any polarization on $X$). Therefore  $Y$ is Seshadri $A$-big in $X$.

\smallskip

iii) Another interesting  application of Corollary \ref{codim111} starts with the following example due to Serre  (see \cite[p. 232, Example 3.2]{Ample}).
Let $B$ be an elliptic curve over $\mathbb C$,
and consider a non-trivial extension of vector bundles
\begin{equation}\label{e:(2.2.15.1)}\begin{CD}
0@>>>\mathscr O_B@>>>E@>{\grb}>>\mathscr O_B@>>>0.\end{CD}\end{equation}
Indeed, the isomorphism classes of such extensions are classified by
$H^1(B,\mathscr O_B)$, which is a one-dimensional vector space because
$B$ is an elliptic curve. Thus non-trivial extensions of type \eqref{e:(2.2.15.1)}
exist. Consider the geometrically ruled surface $Y:=\mathbb P(E)$ associated
to the rank two vector bundle $E$, and denote by $p\colon Y\to B$ the canonical
projection. Let $C$ be the section of $p$ corresponding to the surjection
$\grb$. Since the extension \eqref{e:(2.2.15.1)} is non-trivial one can show
(\cite[p. 232, Example 3.2]{Ample}) that the following properties hold.
\begin{enumerate}
\item[a)] $\int_YC^{\cdot 2}=0$ and, in particular, $Y\setminus C$ is not an affine open subset of $Y$.
\item[b)] The associated analytic space $(X\setminus Y)^{\an}$ is biholomorphically  isomorphic to
$\mathbb C^* \times \mathbb C^*$, and therefore is a Stein manifold.
\end{enumerate}
Set $N:=\mathscr O_Y(C)$ and $V:=Y\setminus C$. From a) and b) it follows that any irreducible curve $C'$ of $Y$, $C'\neq C$, satisfies the following properties: $C'\cap V\neq\varnothing$ and $C'\not\subseteq V$. Such properties imply that for every irreducible curve $C'\neq C$ on $Y$ one has $\int_Yc_1(N)\cdot c_1(\mathscr O_Y(C'))>0$, i.e., $\deg(N\otimes\mathscr O_{C'})>0$ for every irreducible curve $C'$ different from $C$. Recalling also that the normal bundle of $C$ in $Y$ has degree $0$,  we then  proved that hypothesis \eqref{22} of Corollary 
\ref{codim111} is satisfied. At this point, setting $X=\mathbb P(N^{-1}\oplus\mathscr O_Y)$ and taking $A$ any polarization on $X$, we get a polarized threefold $(X,A)$ containing a divisor $Y$  (the zero section of the bundle) with normal bundle $N=\mathscr O_Y(C)$. Then by Corollary 
\ref{codim111} we conclude that $Y$ is Seshadri $A$-big in $X$.
\end{exs*}

\begin{re}{\em  Further examples  of non-ample line bundles $N$ on a projective manifold $Y$ such that $\deg(N\otimes\mathscr O_C)>0 $ for every irreducible curve $C$ in $Y$ can be found in  \cite[Examples 3.13 and  3.14]{BS1}. Starting with such a pair $(Y,N)$ and performing the construction $X=\mathbb P(N^{-1}\oplus\mathscr O_Y)$ as in Examples \ref{Mumford}, we get an effective divisor $Y$ in $X$ which is Seshadri $A$-big with respect to any polarization $A$ of $X$. Moreover,
in these examples  (including Example \ref{Mumford}, iii), in which  $N$ is not nup in the sense of \cite{BS1})  the projective manifold $X=\mathbb P(N^{-1}\oplus\mathscr O_Y)$  contains $Y$ as the ``zero section''  such that $Y$ is a Seshadri $A$-big divisor in $X$ with respect to any polarization $A$ on $X$. Then the ``section at infinity''  (the image $Z$ of the closed immersion  $Y\subset X=\mathbb P(N^{-1}\oplus\mathscr O_Y)$ corresponding to the canonical surjection $N^{-1}\oplus\mathscr O_Y\to\mathscr O_Y$) is a hypersurface of $X$ such that $Y\cap Z=\varnothing$. However, by Theorem \ref{Main} below, it will follow that in all these examples  $Y$ cannot be Seshadri $A$-ample in $X$.}\end{re}

\section{Seshadri $A$-big submanifolds of higher codimension}\label{exampl}
\addtocounter{subsection}{1}\setcounter{theorem}{0}

In this section  we give a series of relevant examples of submanifolds $Y$ of dimension $\geq 2$ and of codimension $\geq 2$ in a polarized manifold $(X,A)$  that are
Seshadri $A$-big, giving a good motivation for the study of this notion in general.   Examples of Seshadri $A$-big curves in a polarized $k$-fold can be found in \cite{Pao}, if $k=3$,   and subsequently in \cite{BBF}, if $k\geq 4$. Proposition  \ref{finite} at the end of this section will show that Seshadri $A$-bigness is stable  under finite coverings.
\begin{prgrph*}{Complete intersections in codimension $2$}\label{CI'}  Assume that $Y$ is the  complete intersection surface in a fourfold $X$, embedded in $\mathbb P^n$, with two hypersurfaces of $\mathbb P^n$ of degrees $d_1\geq d_2$. Let $A=\mathscr O_X(1)$.  By Example \ref{CI}, we have $\varepsilon=\varepsilon(Y,A)=\frac{1}{d_1}$. Moreover, $d:=\int_Yc_1(\mathscr  O_Y(1))^{\cdot2}$ is the degree of $Y$. Then equality \eqref{y=2}  immediately gives
$$
\delta(Y,A)=d(-\varepsilon^2(d_1^2+d_2^2+d_1d_2)+2\varepsilon(d_1+d_2)-1)=d\Big(\frac{d_2}{d_1}-\frac{d_2^2}{d_1^2}\Big).$$ 
Therefore $\delta(Y,A)>0$ if $d_1>d_2$, and $\delta(Y,A)=0$ if $d_1=d_2$. In the former case we get $\delta_{\eta}(Y,A)>0$, and in the latter case we have
$$\delta_{\eta}(Y,A)=d(-3d_1^2\eta^2+4d_1\eta-1).$$
Thus the polynomial $f(x)=-3d_1^2x^2+4d_1x-1$ (with real coefficients), which has the roots $\frac{1}{d_1}$ and $\frac{1}{3d_1}$, assumes positive values for every $x\in\big(\frac{1}{3d_1}, \frac{1}{d_1}\big)$. Then, for every $d_1\geq d_2$, we have $\delta_{\eta}(Y,A)>0$ for every $\eta\in\big(\frac{1}{3d_1}, \frac{1}{d_1}\big)$. 

\smallskip

 Assume now that $Y$ is the  complete intersection threefold of the fivefold $X$, embedded in $\mathbb P^n$, with two hypersurfaces of $\mathbb P^n$ of degrees 
$d_1\geq d_2$. Let $A=\mathscr O_X(1)$.  By Example \ref{CI} again,  we have $\varepsilon(Y,A)=\frac{1}{d_1}$. Let $d:=\int_Yc_1(\mathscr  O_Y(1))^{\cdot3}$ be the degree of $Y$. Proceeding similarly as above and using \eqref{y=3} we immediately get:
$$\delta(Y,A)=d\left(\Big(\frac{d_2}{d_1}\Big)^3-2\Big(\frac{d_2}{d_1}\Big)^2+\frac{d_2}{d_1}\right)=\frac{dd_2}{d_1}\Big(\frac{d_2}{d_1}-1\Big)^2.$$
It follows that $\delta(Y,A)\geq 0$, with equality if and only if $d_1=d_2$. On the other hand,
$$\delta_{\eta}(Y,A)=d\big(-1+3(d_1+d_2)\eta-3(d_1^2+d_1d_2+d_2^2)\eta^2+(d_1^3+d_1^2d_2+d_1d_2^2+d_2^3)\eta^3\big).$$
Then if in the polynomial (with real coefficients)
$$(d_1^3+d_1^2d_2+d_1d_2^2+d_2^3)x^3-3(d_1^2+d_1d_2+d_2^2)x^2+3(d_1+d_2)x-1$$
we put $d_1=d_2$, we get the polynomial
$$P(x)=4d_1^3-9d_1^2x^2+6d_1x-1=4d_1^3\Big(x-\frac{1}{d_1}\Big)^2\Big(x-\frac{1}{4d_1}\Big).$$
As above,  it follows that $\delta_{\eta}(Y,A)>0$ for every $\eta\in(\frac{1}{4d_1},\frac{1}{d_1})$. 

\smallskip

Recalling Lemma \ref{Big},  the previous computations leads to the following result: 
{\em if $Y$ is a smooth complete intersection of codimension $2$ of a projective submanifold  $X$ of dimension $k=4,5$ of $\mathbb P^n$,  then $Y$ is Seshadri $\mathscr O_X(1)$-big}.
However, from a more conceptual point of view, we shall generalize this result to all complete intersections of codimension $2$, see
Proposition \ref{CIcod2} below.

\smallskip

Notice that in codimension $k-y\geq 3$ not every complete intersection is Seshadri 
$\mathscr O_X(1)$-big. This already happens  for curves. Precisely, if $Y$ is a smooth complete intersection curve in $X$ ($X\subseteq\mathbb P^n$) of multidegree $d_1\geq d_2\geq\cdots\geq d_{k-1}\geq 1$, then $Y$ is Seshadri $\mathscr O_X(1)$-big if and only if $d_2+\cdots+d_{k-1}>(k-3)d_1$, see \cite[Example 1.7]{BBF}.

Assume now that $(k,y)=(5,2)$, i.e., $Y$ is a smooth complete intersection surface of the $5$-fold $X\subset\mathbb P^5$ with three hypersurfaces $\mathbb P^5$ of degrees $d_1\geq d_2\geq d_3\geq 1$. Then the Seshadri constant is $\varepsilon(Y,\mathbb P^5)=\frac{1}{d_1}$ and by \eqref{Ses2'},  for every $\eta\in(0,\frac{1}{d_1})\cap\mathbb Q$, we have
$$\delta_{\eta}(Y,\mathscr O_X(1))=d\eta\big(-(d_1^2+d_2^2+d_3^2+d_1d_2+d_1d_3+d_2d_3)\eta^2+3(d_1+d_2+d_3)\eta-3\big).$$
Consider  the polynomial function $f\colon\mathbb R\to\mathbb R$ defined by
\begin{equation}\label{pol2} f(x)=-(d_1^2+d_2^2+d_3^2+d_1d_2+d_1d_3+d_2d_3)x^2+3(d_1+d_2+d_3)x-3. \end{equation}
For example, if we take $X=\mathbb P^5$, $d_1=3$ and $d_2=d_3=2$ we get $f(x)=-33x^2+21x-3$. Then $f(\frac{1}{3})=\frac{1}{3}$, whence $Y$ is Seshadri $\mathscr O_{\mathbb P^5}(1)$-big.

If instead we take $X=\mathbb P^5$, $d_1=9$ and $d_2=d_3=2$ we get $f(x)=-129x^2+39x-3$, whose discriminant is $\Delta=-27<0$. We deduce that $f(x)<0$ for every $x\in\mathbb R$, which implies that in this case $Y$ is never Seshadri $\mathscr O_{\mathbb P^5}(1)$-big.
\end{prgrph*}
\begin{prgrph*}{Zero loci of sections of ample vector bundles}\label{Zero} Fix two integers $k$ and $y$ such that $1\leq y\leq k-2$. Let $\se$ be an ample vector bundle of rank 
 $k-y$ on a smooth projective $k$-dimensional variety $X$, $k\geq 3$, and set
 $A:=\det(\se)$.  As $\se$ is ample, a general result of Hartshorne  \cite{Ample} implies that $A$ is an ample line bundle on $X$. Let $s\in H^0(X,\se)$ be a section
 such that its zero locus $Y:=Z(s)$ is a smooth $y$-dimensional subvariety of $X$. Then by a result
of Sommese, $Y$ is  connected (see \cite[(1.16)]{So}, or also
 \cite{FL}). Then $s$ defines a map $s:\so_X\to \se$ and,
 taking the dual map $s^*:\se^*\to\so_X$, one gets the
 surjection $s^*:\se^*\to\sy_Y$, where $\sy_Y$ is the ideal sheaf of $Y$ in
 $X$. Thus we get a surjective map ${\bf S}(\se^*)\to\oplus_{i=0}^{\infty}\sy_Y^i$,
 where ${\bf S}(\se^*)$ is the symmetric $\so_X$-algebra of $\se^*$. It follows that
 there is a closed immersion $i:X_Y\hookrightarrow {\mathbb P}(\se^*)$ of
 $X$-schemes such that $i^*(\so_{{\mathbb P}(\se^*)}(1))\cong\so_{X_Y}(-E)$.
 Consider the canonical isomorphisms
 $$\wedge^{k-y-1}\se\cong\se^*\otimes\det(\se)=\se^*\otimes A.$$
  It follows that there is
 a canonical isomorphism ${\mathbb P}(\wedge^{k-y-1}\se)\cong{\mathbb P}(\se^*)$ with
 the property that $\so_{{\mathbb P}(\wedge^{k-y-1}\se)}(1)\cong
 \so_{{\mathbb P}(\se^*)}(1)\otimes p^*(A)$, where 
 $p:{\mathbb P}(\se^*)\to X$ is the canonical
 projection. As $\se$ is ample by hypothesis, a result of Hartshorne
 \cite{AVB}, or also \cite{Ample}, shows that $\wedge^{k-y-1}\se$ is also ample, i.e.,
 $\so_{{\mathbb P}(\wedge^{k-y-1}\se)}(1)\cong\so_{{\mathbb
 P}(\se^*)}(1)\otimes p^*(A)$ is ample. Recalling the above
 immersion, this last fact implies that $A^*-E$ is ample on $X_Y$, whence $\gre(Y,A)>1$. 
 
 \smallskip

 Consider first the case when the vector bundle $\se$ is decomposable of rank $2$, i.e., $\se=L'\oplus M'$, with $L'$ and $M'$ ample line bundles on $X$. By what we have said above, $\varepsilon(Y,A)>1$, with $A=\det(\se)=L'\otimes M'$. Moreover,  the normal bundle $N$ is isomorphic to $L\oplus M$, where $L=L'|Y$ and $M=M'|Y$ (in particular, $L$ and $M$ are ample line bundles on $Y$). Then it is easy to compute the Segre classes of the normal bundle $N=\se_Y$. Namely, for every line bundle $\sL$ on $Y$ and for every $i\geq 0$ one has $s_i(\sL)=(-1)^ic_1(\sL)^{\cdot i}$. Therefore for every integer $t\geq 0$ we get
\begin{equation}\label{genoa}s_t(N)=s_t(L\oplus M)=\sum_{i=0}^ts_i(L)\cdot s_{t-i}(M)= (-1)^t\sum_{i=0}^tc_1(L)^{\cdot i}\cdot c_1(M)^{\cdot (t-i)}.\end{equation}
Using \eqref{delta1} and \eqref{genoa} and taking into account that $y-k+t+2=t$ (because 
$\codim_X(Y)=k-y=2$) in our situation we obtain
\begin{eqnarray*}\delta_{\eta}(Y,A)&=&-\sum_{t=0}^{k-2}\binom{k-2}{t}\eta^t\int_Ys_t(N)\cdot
c_1(A_Y)^{\cdot(k-t-2)}\end{eqnarray*}
\begin{eqnarray*}&=&\sum_{t=0}^{k-2}(-1)^{t+1}\binom{k-2}{t}\eta^t\int_Y\sum_{i=0}^t
c_1(L)^{\cdot i}\cdot c_1(M)^{\cdot(t-i)}\cdot c_1(L\otimes M)^{\cdot(k-t-2)}\\
&=&\sum_{t=0}^{k-2}(-1)^{t+1}\binom{k-2}{t}\eta^t\int_Y(c_1(L)+c_1(M))^{\cdot(k-t-2)}\cdot\sum_{i=0}^t
c_1(L)^{\cdot i}\cdot c_1(M)^{\cdot(t-i)}.\\ \label{alg1}
\end{eqnarray*}

To proceed further we need the following elementary result of algebra.

\begin{lemma}\label{algebra} Let $B$ be a commutative unitary ring and let $b_1,b_2\in B$.
Then for every integer $y\geq 1$ one has
$$ \sum_{t=0}^{y}(-1)^{t+1}\binom{y}{t}(b_1+b_2)^{y-t}\sum_{i=0}^t
b_1^{i}b_2^{t-i}= \sum_{j=1}^{y-1}b_1^{y-j}b_2^j.$$
\end{lemma}

\proof Consider the polynomial ring $\mathbb Z[T_1,T_2]$ in two variables $T_1$ and $T_2$ over the ring of integers $\mathbb Z$. We claim that the identity 
\begin{equation}\label{pol}\sum_{t=0}^{y}(-1)^{t+1}\binom{y}{t}(T_1+T_2)^{y-t}\sum_{i=0}^t
T_1^{i}T_2^{t-i}= \sum_{j=1}^{y-1}T_1^{y-j}T_2^j\end{equation}
holds true. 
Indeed, if we denote by $I$ the left hand side member of the identity to be proved, we have successively
(in the fraction field $\mathbb Q(T_1,T_2)$):
\begin{eqnarray*}
I&=&\sum_{t=0}^{y}(-1)^{t+1}\binom{y}{t}(T_1+T_2)^{y-t}\;\frac{T_1^{t+1}-T_2^{t+1}}{T_1-T_2}\\&=&
\frac{1}{T_1-T_2}\sum_{t=0}^{y}(-1)^{t+1}\binom{y}{t}(T_1+T_2)^{y-t}T_1^{t+1}-\frac{1}{T_1-T_2}\sum_{t=0}^{y}(-1)^{t+1}\binom{y}{t}(T_1+T_2)^{y-t}T_2^{t+1}\\
&=&\frac{-T_1}{T_1-T_2}\sum_{t=0}^{y}\binom{y}{t}(T_1+T_2)^{y-t}(-T_1)^{t}+\frac{T_2}{T_1-T_2}\sum_{t=0}^{y}\binom{y}{t}(T_1+T_2)^{y-t}(-T_2)^{t}\\
&=&\frac{-T_1}{T_1-T_2}\big((T_1+T_2)-T_1\big)^y+\frac{T_2}{T_1-T_2}\big((T_1+T_2)-T_2\big)^y\\
&=&\frac{T_1T_2(T_1^{y-1}-T_2^{y-1})}{T_1-T_2}=T_1T_2\sum_{j=1}^{y-1}T_1^{y-j-1}T_2^{j-1}
=\sum_{j=1}^{y-1}T_1^{y-j}T_2^j,
\end{eqnarray*}
which proves our claim. 

To prove the identity in  general, consider the ring homomorphism
$$\varphi\colon\mathbb Z[T_1,T_2]\to B$$ 
such that $\varphi(n)=n\cdot 1_B$ and $\varphi(T_i)=b_i$, $i=1,2$. Then the conclusion follows from \eqref{pol} because 
\begin{multline*} \sum_{t=0}^{y}(-1)^{t+1}\binom{y}{t}(b_1+b_2)^{y-t}\sum_{i=0}^tb_1^{i}b_2^{t-i}-\sum_{j=1}^{y-1}b_1^{y-j}b_2^j\\
=\varphi\Big(\sum_{t=0}^{y}(-1)^{t+1}\binom{y}{t}(T_1+T_2)^{y-t}\sum_{i=0}^t
T_1^{i}T_2^{t-i} -\sum_{j=1}^{y-1}T_1^{y-j}T_2^j\Big)=\varphi(0)=0.\end{multline*}\vskip -1cm\qed

\vskip 0.3cm

Now we are ready to prove the following result.

\begin{prop}\label{CIcod2} Let $\se=L'\oplus M'$ be a decomposable ample vector bundle of rank $2$ on a smooth projective variety $X$ of dimension $k\geq 3$, and let $s=(f,g)\in H^0(\se)=H^0(L'\oplus M')$ be a global section such that the zero locus $Y:=Z(s)$ is smooth of codimension $2$ in $X$. Then $Y$ is Seshadri $A$-big, with $A:=L'\otimes M'$. In particular, if we fix a  projective embedding $X\hookrightarrow\mathbb P^n$, every smooth complete intersection in $X$  of codimension $2$ is Seshadri $\mathscr O_X(1)$-big.
\end{prop}

\proof Since we have seen above that the Seshadri constant $\varepsilon(Y,L'\otimes M')$ is $>1$, it will be enough to show that $\delta_1(Y,A)>0$, where $A=L'\otimes M'$. To this end, in Lemma \ref{algebra} take as $B$ the Chow ring $A(Y)$ of cycles of $Y$ modulo the rational equivalence, $y=\dim(Y)$, and put $b_1=c_1(L)$ and $b_2=c_1(M)$ (recall that $L=L'|Y$ and $M=M'|Y$). Then by Lemma \ref{algebra} we get
$$ \Big(\sum_{t=0}^{y}(-1)^{t+1}\binom{y}{t}(c_1(L)+c_1(M))^{\cdot(y-t)}\Big)\cdot \Big(\sum_{i=0}^{t}
c_1(L)^{\cdot i}\cdot c_1(M)^{\cdot(t-i)}\Big)= \sum_{j=1}^{y-1}c_1(L)^{\cdot(y-j)}\cdot c_1(M)^{\cdot j},$$
which implies
$$\delta_1(Y,A)=\sum_{j=1}^{y-1}\int_Yc_1(L)^{\cdot(y-j)}\cdot c_1(M)^{\cdot j}>0,$$
because $\int_Yc_1(L)^{\cdot(y-j)}\cdot c_1(M)^{\cdot j}>0$ for every $j=1,\ldots,y-1$ (since $L$ and $M$ are ample line bundles on $Y$).

To show the last statement take $L'=\mathscr O_X(d_1)$ and $M'=\mathscr O_X(d_2)$, with $d_1$ and $d_2$ positive integers, and use the argument above together with Remark \ref{indep}. \qed

\smallskip

Now let us consider the case when the vector bundle $\se$ is indecomposable. As above, we have $y=k-2$ and $\varepsilon(Y,A)> 1$,  so can take $\eta=1$.
Taking into account of relations \eqref{delta4}, \eqref{y=2}, \eqref{y=3} and \eqref{y=4} and using the facts that $N=\se_Y=\se\otimes\mathscr O_Y$ and $A=\det(\se)$ (whence $A_Y=\det(\se_Y)$ and $c_1(\se_Y)=c_1(A_Y)$), we find:
\begin{eqnarray} \delta(Y,A)&=&\varepsilon(Y,A)\deg(\se_Y)-\deg(\se_Y)>0,\;\;\text{for $(k,y)=(3,1)$,}\label{delta20}\\
\delta_1(Y,A)&=&\int_Yc_2(\se_Y),\;\;\text{for $(k,y)=(4,2)$},\label{delta8}\\
\delta_1(Y,A)&=&\int_Yc_1(\se_Y)\cdot c_2(\se_Y),\;\;\text{for $(k,y)=(5,3)$},\label{delta9}\\
\delta_1(Y,A)&=&\int_Yc_2(\se_Y)\cdot(c_1(\se_Y)^{\cdot2}- c_2(\se_Y)),\;\;\text{for $(k,y)=(6,4)$},\label{delta10}\\
\delta_1(Y,A)&=&\int_Yc_1(\se_Y)\cdot c_2(\se_Y)\cdot (c_1(\se_Y)^{\cdot 2}- 2c_2(\se_Y)),\;\;\text{for $(k,y)=(7,5)$}.\label{delta11}
\end{eqnarray}

Now, as the vector bundle $\se$ is ample, the restriction $\se_Y$ is also ample. We claim that $\delta_1(Y,A)>0$ if $k\leq 6$.  Indeed, if $(k,y)=(4,2)$, this follows from \eqref{delta8} and from a result of Kleiman \cite {Kl} (see also Bloch and Gieseker  \cite{BG} for a subsequent  more general result) according to which $\int_Yc_2(\se_Y)>0$ when $Y$ is a surface. 

If $(k,y)=(5,3)$ by \eqref{delta9} we have
\begin{eqnarray*}\delta_1(Y,A)&=&\int_Yc_1(\se_Y)\cdot c_2(\se_Y)=\int_Y\left|\begin{array}{ccc}
c_2(\se_Y)&0&0\\
1&c_1(\se_Y)&c_2(\se_Y)\\
0&0&1\\
\end{array}\right|
\\&=&
\int_Y\left|\begin{array}{ccc}
c_2&c_3&c_4\\
c_0&c_1&c_2\\
c_{-2}&c_{-1}&c_0\\
\end{array}\right|(\se_Y),\end{eqnarray*}
where the latter determinant is, in terminology of \cite{FL}, the Schur polynomial (evaluated at $\se_Y$) associated  to the partition $r\geq\lambda_1\geq\lambda_2\geq\lambda_3\geq 0$, whith $r=\rank(\se_Y)=2$ and  $(\lambda_1,\lambda_2,\lambda_3)=(2,1,0)$. This  determinant  is strictly positive by the main result of \cite{FL}, or also  \cite[Theorem 8.3.9]{Laz1}.

If $(k,y)=(6,4)$ by \eqref{delta10} we have
\begin{eqnarray*}\delta_1(Y,A)&=&\int_Yc_2(\se_Y)\cdot(c_1(\se_Y)^{\cdot2}- c_2(\se_Y))\\
&=&\int_Y\left|\begin{array}{cccc}
c_2(\se_Y)&0&0&0\\
1&c_1(\se_Y)&c_2(\se_Y)&0\\
0&1&c_1(\se_Y)&c_2(\se_Y)\\
0&0&0&1\\
\end{array}\right|=
\int_Y\left|\begin{array}{cccc}
c_2&c_3&c_4&c_5\\
c_0&c_1&c_2&c_3\\
c_{-1}&c_0&c_1&c_2\\
c_{-3}&c_{-2}&c_{-1}&c_0\\
\end{array}\right|(\se_Y),\end{eqnarray*}
where the latter determinant is, in terminology of \cite{FL}, the Schur polynomial (evaluated at $\se_Y$) associated  to the partition $r\geq\lambda_1\geq\lambda_2\geq\lambda_3\geq\lambda_4\geq 0$, with $r=\rank(\se_Y)=2$ and  $(\lambda_1,\lambda_2,\lambda_3,\lambda_4)=(2,1,1,0)$. Again, the determinant is  strictly positive because $\se_Y$ is ample, by the main result of \cite{FL}, see also \cite[Theorem 8.3.9]{Laz1}.

Assume now  $(k,y)=(7,5)$. The polynomial $c_1\cdot c_2\cdot (c_1^{\cdot 2}- 2c_2)$ occurring in \eqref{delta11} is not positive (because $c_1^2-2c_2$ is not a linear combination of the Schur polynomials $c_1^2-c_2$ and $c_2$, see \cite{FL}), whence  $\delta_1(Y,A)$ is not always positive. Note incidently that the polynomial $c_1^{\cdot 2}- 2c_2$ has a rather interesting story, see \cite[Example 8.3.11, p. 121]{Laz1}, II, and also \cite[p. 418]{GH}. 

Putting things together  and using  Lemma \ref{Big} together with formula \eqref{delta11}, we proved the following result (the case $(k,y)=(3,1)$ of Proposition \ref{Z(s)} was treated in \cite{Pao}).

\begin{prop}\label{Z(s)} Let $\se$ be an ample vector bundle of rank $2$ on a projective manifold $X$ of dimension $k$, with $3\leq k\leq 6$, and let $s\in H^0(X,\se)$ be a global section of $\se$ such that the zero locus $Y:=Z(s)$ is a smooth $2$-codimensional submanifold of $X$. Then $Y$ is Seshadri $A$-big in $X$, where $A:=\det(\se)$. If $(k,y)=(7,5)$ then $Y$ is Seshadri $A$-big if $\int_Yc_1(\se_Y)\cdot c_2(\se_Y)\cdot (c_1(\se_Y)^{\cdot 2}- 2c_2(\se_Y))>0$.\end{prop}
\end{prgrph*}
\begin{prgrph*}{Submanifolds of $\mathbb P^k$ and their Chern classes}\label{Sub} In order to give some  relevant examples of Seshadri $\mathscr O_{\mathbb P^k}(1)$-big submanifolds $Y$ of  $X=\mathbb P^k$ of dimension $y\geq 2$,  we need some generalities regarding Chern classes which allow us to explicitly compute  Chern classes of the normal bundle of some submanifolds  of $\pn k$.

The restriction to $Y$ of the Euler sequence of $\mathbb P^k$,
$$0\to \mathscr O_Y\to\mathscr O_Y(1)^{\oplus k+1}\to T_{\mathbb P^k}|Y\to 0,$$
where $T_{\mathbb P^k}$ is the tangent bundle of $\mathbb P^k$, yields the equality of Chern polynomials (see \cite[p. 50]{Fu})
\begin{equation}\label{chern1}c_t(T_{\mathbb P^k}|Y)=c_t(\mathscr O_Y(1)^{\oplus k+1}).\end{equation}
By additivity,  formula \eqref{chern1} can be rewritten as
$$c_t(T_{\mathbb P^k}|Y)=(1+tc_1(\mathscr O_Y(1))^{k+1}.$$
Taking Chern classes in the latter formula we get
\begin{equation}\label{chern2}c_i(T_{\mathbb P^k}|Y)=\binom{k+1}{i}c_1(\mathscr O_Y(1))^{\cdot i},\;i=0,1,\ldots,y. \end{equation}
In particular, if $i=y$ we get $c_y(T_{\mathbb P^k}|Y)=\binom{k+1}{y}\deg(Y)$.

On the other hand, the normal sequence 
$$0\to T_Y\to T_{\mathbb P^k}|Y\to N\to 0$$ and the additivity of the Chern classes yield the equality
$$c_t(T_{\mathbb P^k}|Y)=c_t(T_Y)\cdot c_t(N),$$
where $N$ is the normal bundle of $Y$ in ${\mathbb P^k}$. The latter equality can be rewritten via \eqref{chern2} as
\begin{equation}\label{chern3}\sum_{i=0}^y\binom{k+1}{i}c_1(\mathscr O_Y(1))^{\cdot i}t^i=\Big(\sum_{i=0}^yc_i(Y)t^i\Big)\cdot \Big(\sum_{i=0}^yc_i(N)t^i\Big),\end{equation}
where $c_i(Y):=c_i(T_Y)$, $i=0,1,\ldots,y$, are the Chern classes of $Y$.
\end{prgrph*}
\begin{prgrph*}{Rational normal scrolls of dimension $2$}\label{RatScroll} Let $Y_e$ be the image in $\mathbb P^{e+3}$ of the geometrically  ruled surface (or Segre--Hirzebruch surface) $\mathbb F_e:=\mathbb P(\mathscr O_{\mathbb P^1}\oplus\mathscr O_{\mathbb P^1}(-e))$, $e\geq 0$, embedded as a scroll in $\mathbb P^{e+3}$ via the complete linear system $|C_0+(e+1)F|$, where  $C_0$ is the minimal section (with $C_0^{\cdot 2}=-e$) and $F$ is any fiber of the canonical projection $\pi\colon\mathbb F_e\to\mathbb P^1$ respectively (we follow here the notation as in  Hartshorne \cite{Hartshorne}). 
Moreover, $Y_e$ is a surface of minimal degree $e+2$ in $\mathbb P^{e+3}$, whence by Example \ref{MinDeg} we have
\begin{equation}\label{min}\varepsilon(Y_e,\mathscr O_{\mathbb P^{e+3}}(1))=\frac{1}{2}.\end{equation}

Furthermore,  $c_1(Y_e)= -K_e$ and $\int_Yc_2(Y_e)=\chi_{\top}(Y_e)=2-2b_1+b_2=4$ (see \cite[Chapter 5]{B1}), where $K_e$ is the canonical class of $Y_e$, 
$\chi_{\top}(Y_e)$ is the topological Euler characteristic of $Y_e$, and $b_i$, $i=1,2$, are the Betti numbers of $Y_e$. In our case, $b_1=0$ and $b_2=2$, so $\int_{Y_e}c_2(Y_e)=4$. 

Let $N$ be the normal bundle of $Y_e$ in $\mathbb P^{e+3}$. Then  the identity \eqref{chern3} gives
$$c_1(\mathscr O_{Y_e}(-K_e))+c_1(N)=(e+4)c_1(\mathscr O_{Y_e}(1))=c_1(\mathscr O_{Y_e}(e+4)),$$
 and
$$\int_{Y_e}c_2(N)-\int_{Y_e}c_1(\mathscr O_{Y_e}(K_e))\cdot c_1(N)+4=\binom{e+4}{2}\deg(Y_e)=\frac{(e+4)(e+3)(e+2)}{2}.$$
 Recalling that $\mathscr O_{Y_e}(K_e)=\mathscr O_{Y_e}(-2C_0-(e+2)F)$ and 
$\mathscr O_{Y_e}(e+4)=\mathscr O_{Y_e}((e+4)C_0+(e+4)(e+1)F)$, from the first equality we get 
\begin{equation}\label{chern4}c_1(N)=c_1(\mathscr O_{Y_e}((e+2)C_0+(e^2+4e+2)F)),\end{equation}
while the second equality yields
\begin{equation}\label{chern5}\int_{Y_e}c_2(N)=\frac{(e+4)(e+3)(e+2)}{2}-e^2-8e-12.\end{equation}
  Now using \eqref{chern4} and \eqref{chern5}, we get
$$\int_{Y_e}(c_1(N)^{\cdot2}-c_2(N))=\frac{e^3+9e^2+22e+16}{2}\;\;\text{and}\;\;
\int_{Y_e}c_1(N)\cdot c_1(\mathscr O_{Y_e}(1))=e^2+5e+4.$$
Therefore taking into account of \eqref{min} and \eqref{Ses2} (with $k=e+3$ and $\deg(Y_e)=e+2$) we find, for every $e\geq 0$,
\begin{equation}\label{Ses5}\delta_{\eta}(Y_e,\mathscr O_{\mathbb P^{e+3}}(1))=\frac{\eta^{e-1}}{2}\big(-(e^3+9e^2+22e+16)\eta^2+2(e+1)(e^2+5e+4)\eta-e(e+1)(e+2)\big).\end{equation}
Letting $\eta=\varepsilon(Y_e,\mathscr O_{\mathbb P^{e+3}}(1))=\frac{1}{2}$, after a straightforward  computation the latter equality yields:
\begin{equation}\label{Ses6}\delta(Y_e,\mathscr O_{\mathbb P^{e+3}}(1))=\frac{e}{2^{e+2}}(-e^2+3e+6).\end{equation}

We can now prove  the following:

\begin{prop}\label{F1} The image $Y_e$ in $\pn {e+3}$ of the rational normal scroll $\mathbb F_e$   is Seshadri $\mathscr O_{\mathbb P^{e+3}}(1)$-big if  and only if $0\leq e\leq 4$. 
\end{prop}

\proof From \eqref{Ses6} it follows that $\delta(Y_e,\mathscr O_{\mathbb P^{e+3}}(1))>0$ if and only if $0\leq e\leq 4$. It follows that $Y_e$ is Seshadri $\mathscr O_{\mathbb P^{e+3}}(1)$-big if $0\leq e\leq 4$.

On the other hand,  racalling \eqref{Ses5}, consider the function $f\colon\mathbb R\to\mathbb R$  defined by
$$f(x)=-(e^3+9e^2+22e+16)x^2+2(e+1)(e^2+5e+4)x-e(e+1)(e+2).$$
Then $f$ is increasing on the interval $(0,\varepsilon(Y_e,\mathscr O_{\mathbb P^{e+3}}(1)))=(0,\frac{1}{2})$ if $e\geq 2$.
Indeed, the function $f$ assumes the maximum value on $\mathbb R$  for 
$$x=\frac{(e+1)(e^2+5e+4)}{e^3+9e^2+22e+16}=\frac{e^3+6e^2+9e+4}{e^3+9e^2+22e+16},$$
and, moreover (as it is easily checked), 
$$\frac{e^3+6e^2+9e+4}{e^3+9e^2+22e+16}>\frac{1}{2}\;\;\text{ if $e\geq 2$.}$$
This implies that  for $e\geq 2$ the function $f$ is increasing on the interval $(0,\frac{1}{2})$. It follows that if $e\geq 2$ the surface $Y_e$ is Seshadri $\mathscr O_{\mathbb P^{e+3}}(1)$-big in $\mathbb P^{e+3}$ if  and only if $\delta(Y_e,\mathscr O_{\mathbb P^{e+3}}(1))>0$, or else (by the remark made  at the beginning), only if $e\leq 4$. \qed
\end{prgrph*}
\begin{prgrph*}{The Veronese surface in $\mathbb P^5$ and its projection in $\mathbb P^4$}\label{Vero} The last example of smooth surface of minimal degree is the Veronese surface $Y$ in $\mathbb P^5$, i.e., the image of the Veronese embedding $v_2\colon\mathbb P^2\hookrightarrow\mathbb P^5$. By Example \ref{MinDeg}, the Seshadri constant is $\varepsilon(Y,\mathscr O_{\mathbb P^5}(1))=\frac{1}{2}$ whence, by formula 
\eqref{Ses2},
\begin{eqnarray*}\delta(Y,\mathscr O_{\mathbb P^5}(1))&=&-\frac{1}{8}\int_Y(c_1(N)^{\cdot 2}-c_2(N))+\frac{3}{4}\int_Yc_1(N)\cdot c_1(\mathscr O_Y(1))-\frac{3}{2}\deg(Y)\\
&=&-\frac{1}{8}\int_Y(c_1(N)^{\cdot 2}-c_2(N))+\frac{3}{4}\int_Yc_1(N)\cdot c_1(\mathscr O_{\mathbb P^2}(2))-6.
\end{eqnarray*}
Computing the Chern classes of the normal bundle $N$   by using \eqref{chern3} (as in the case of the surface scrolls),  we immediately find
$$c_1(N)=c_1(\mathscr O_{\mathbb P^2}(9)) \;\;\text{and}\;\;\int_Yc_2(N)=30.$$
Putting things together we get
\begin{equation}\label{ver}\delta(Y,\mathscr O_{\mathbb P^5}(1))=-\frac{51}{8}+\frac{27}{2}-6=\frac{9}{8}>0.\end{equation}

Consider now the linear projection in $\mathbb P^4$ of the Veronese surface $Y=v_2(\mathbb P^2)$ from a general point of $\mathbb P^5$. It is a surface $Y'\cong\mathbb P^2$ in $\mathbb P^4$ such that $\mathscr O_{Y'}(1)=\mathscr O_{\mathbb P^2}(2)$. Let $N':=N_{Y'|\pn 4}$ be the normal bundle of $Y'$ in $\pn 4$. Then it is known that the homogeneous ideal $\mathscr I_+(Y')$ of $Y'$ in $\mathbb P^4$ is generated by homogeneous polynomials of degree $3$. (We thank John Abbott for providing us an explicit set of generators of degree $3$ of $\mathscr I_+(Y')$ using CoCoA-4.7.) Then by Corollary 1.6 in \cite{BBF}, $\varepsilon(Y',\mathscr O_{\mathbb P^4}(1))\geq \frac{1}{3}$.
Therefore for $\eta\in(0,\frac{1}{3}]$ we have 
\begin{eqnarray*}\delta_{\eta}(Y',\mathscr O_{\mathbb P^4}(1))&=&-\eta^2\int_{Y'}(c_1(N')^{\cdot2}-c_2(N'))+2\eta\int_{Y'}c_1(N')\cdot c_1(\mathscr O_{Y'}(1))-\deg(Y')\\&=&-\eta^2\int_{Y'}(c_1(N')^{\cdot2}-c_2(N'))+2\eta\int_{Y'}c_1(N')\cdot c_1(\mathscr O_{\mathbb P^2}(2))-4.\end{eqnarray*}
Computing as above the Chern classes of $N'$ we find
$$c_1(N')=c_1(\mathscr O_{\mathbb P^2}(7))\;\;\text{and}\;\;\int_{Y'}c_2(N')=16,$$
whence
$$\delta_{\eta}(Y',\mathscr O_{\mathbb P^4}(1))=-33\eta^2+28\eta-4.$$
In particular,

\begin{equation}\label{pr}\delta_{\frac{1}{3}}(Y',\mathscr O_{\mathbb P^4}(1))=-\frac{33}{9}+\frac{28}{3}-4=\frac{5}{3}>0.\end{equation}

Summarizing, \eqref{ver} and \eqref{pr} yield the following result.

\begin{prop}\label{Veronese2} The Veronese surface $($i.e., the image $Y$ of the Veronese embedding $\mathbb P^2\hookrightarrow\mathbb P^5)$ is Seshadri $\mathscr O_{\mathbb P^5}(1)$-big
in $\mathbb P^5$. Let $Y'\cong\mathbb P^2$ be the image of $Y$ via the linear projection in $\mathbb P^4$ of $Y$ from a general point of $\mathbb P^5$. Then $Y'$ is Seshadri $\mathscr O_{\mathbb P^4}(1)$-big in $\mathbb P^4$.
\end{prop}
\end{prgrph*}
\begin{prgrph*}{Some Segre embeddings}\label{SegreEmb} We prove here the $\mathscr O(1)$-bigness of some Segre embeddings.

\smallskip
First, consider the Segre embedding  $i_5: \mathbb P^2\times\mathbb P^1\hookrightarrow\mathbb P^5$, and 
set $Y:=i_5(\mathbb P^2\times\mathbb P^1)$. By Example \ref{MinDeg}, $\varepsilon(Y,\mathscr O_{\mathbb P^5}(1))=\frac{1}{2}$. Since 
$\codim_{\mathbb P^5}(Y)=2$ we can apply  \eqref{y=3}  to get
\begin{multline}\label{scrolldim31}\delta_{\eta}(Y,\mathscr O_{\mathbb P^5}(1))=-3+3\eta\int_Yc_1(N)\cdot c_1(\mathscr O_Y(1))^{\cdot2}\\ -3\eta^2(c_1(N)^{\cdot2}- c_2(N))\cdot c_1(\mathscr O_Y(1))-\eta^3\int_Y(2c_1(N)\cdot c_2(N)-c_1(N)^{\cdot3}),\end{multline}
where $\delta_{\eta}:=\delta_{\eta}(Y,\mathscr O_{\mathbb P^5}(1))$. The equation \eqref{chern3} reads in this case
\begin{equation}\label{chern3'}(1+c_1(\mathscr O_Y(1))t)^6=\big(1+c_1(Y)t+c_2(Y)t^2+c_3(Y)t^3\big)\cdot\big(1+c_1(N)t+c_2(N)t^2\big).\end{equation}
Identifying the coefficients of $t$ and using $c_1(Y)=-c_1(K_Y)$, we get
$$c_1(N)=c_1(K_Y)+6c_1(\mathscr O_Y(1))=c_1(\mathscr O_Y(6))+c_1(K_Y),$$
whence
\begin{equation}\label{normal1}c_1(N)=c_1(\mathscr O(3,4)),\end{equation}
taking into account that $\mathscr O_Y(1)=\mathscr O(1,1)$ and $K_Y=\mathscr O(-3,-2)$, where as usual we set $\mathscr O(a,b):=p_1^*(\mathscr O_{\mathbb P^2}(a))\otimes 
p_2^*(\mathscr O_{\mathbb P^1}(b))$, with $p_1$ and $p_2$ the canonical projections of $\mathbb P^2\times\mathbb P^1$.
Identifying the coefficients of $t^2$ in \eqref{chern3'} we also get
\begin{equation}\label{n2}c_2(N)-c_1(K_Y)\cdot c_1(N)+c_2(T_Y)=15 c_1(\mathscr O_Y(1))^{\cdot 2}. \end{equation}
Using \eqref{normal1} and the fact that 
$$T_Y=p_1^*(T_{\mathbb P^2})\oplus p_2^*(T_{\mathbb P^1})$$ we find
$$c_t(T_Y))=c_t(p_1^*(T_{\mathbb P^2}))\cdot c_t(p_2^*(T_{\mathbb P^1}))=p_1^*(c_t(T_{\mathbb P^2}))\cdot p_2^*(c_t(T_{\mathbb P^1})).$$
Since $c_t(T_{\mathbb P^2})=1+c_1(\mathscr O_{\mathbb P^2}(3))t+ 3xt^2$, with $x\in\mathbb P^2$ and $c_t(T_{\mathbb P^1})=1+c_1(\mathscr O_{\mathbb P^1}(2))t$, we get
$$c_t(T_Y)=(1+c_1(\mathscr O(3,0))t+3(x\times\mathbb P^1)t^2)\cdot(1+c_1(\mathscr O(0,2)t).$$
Identifying the coefficients of $t^2$ in this latter identity yields
\begin{equation}\label{n3} c_2(T_Y)=3(x\times\mathbb P^1)+c_1(\mathscr O(3,0))\cdot c_1(\mathscr O(0,2))=3(x\times\mathbb P^1)+6(\ell\times y),\end{equation}
where $\ell$ is a line of $\mathbb P^2$  and $y\in\mathbb P^1$. 
Substituting   \eqref{n3} in \eqref{n2} and using the obvious formula 
\begin{equation}\label{error}c_1(\mathscr O(a,b))\cdot c_1(\mathscr O(a',b'))=aa'(x\times\mathbb P^1)+(ab'+a'b)(\ell\times y),\end{equation}
we get
\begin{equation}\label{n5} c_2(N)=3 (x\times\mathbb P^1)+6 (\ell\times y),\;\;\text{with $x\in\mathbb P^2$, $y\in\mathbb P^1$, and $\ell$ a line in $\mathbb P^2$}. 
\end{equation}

Now we are ready to compute $\delta_{\eta}(Y,\mathscr O_{\mathbb P^5}(1))$. Using repeatedly \eqref{normal1}, \eqref{n5} and \eqref{error} we have:
\begin{eqnarray*}
\int_Yc_1(N)\cdot c_1(\mathscr O_Y(1))^{\cdot2}  &=&\int_Yc_1(\mathscr O(3,4))\cdot c_1(\mathscr O(1,1))\cdot c_1(\mathscr O(1,1)) \\& =& 
\int_Yc_1(\mathscr O(3,4))\cdot[(x\times\mathbb P^1)+2(\ell\times y)]=10.\\
\int_Y(c_1(N)^{\cdot2}-c_2(N))\cdot c_1(\mathscr O_Y(1))&=&\int_Yc_1(\mathscr O(3,4))\cdot c_1(\mathscr O(3,4))\cdot c_1(\mathscr O(1,1)) \\ & &
-\int_Y[3(x\times\mathbb P^1)+6(\ell\times y)]\cdot c_1(\mathscr O(1,1))=24.\\
\int_Y(2c_1(N)\cdot c_2(N)-c_1(N)^{\cdot3})&=&2\int_Yc_1(\mathscr O(3,4))\cdot[3(x\times\mathbb P^1)+6(\ell\times y)] \\ & &
-\int_Yc_1(\mathscr O(3,4))\cdot[9(x\times\mathbb P^1)+24(\ell\times y)] =-48.
\end{eqnarray*}
Substituting these equalities in \eqref{scrolldim31} we find
$$\delta_{\eta}(Y,\mathscr O_{\mathbb P^5}(1))=48\eta^3-72\eta^2+30\eta-3.$$
Taking $\eta=\frac{1}{3}\in(0,\frac{1}{2})$ we get
$$\delta_{\frac{1}{3}}(Y,\mathscr O_{\mathbb P^5}(1))=\frac{48}{27}-\frac{72}{9}+\frac{30}{3}-3=\frac{16}{9}-8+10-3=\frac{16}{9}-1=\frac{7}{9}>0.$$
Note that $\delta(Y,\mathscr O_{\mathbb P^5}(1))=\delta_{\frac{1}{2}}(Y,\mathscr O_{\mathbb P^5}(1))=0$. (Using this and the fact that $f'(\frac{1}{2})<0$ we also deduce that $f(\eta)>0$ for every $\eta<\frac{1}{2}$ which is close to $\frac{1}{2}$, where $f(\eta)=48\eta^3-72\eta^2+30\eta-3$.)
Therefore $Y:=i_5(\mathbb P^2\times\mathbb P^1)$ is $\mathscr O_{\mathbb P^5}(1)$-big.

\smallskip

Next, consider the Segre embedding  $i_7: \mathbb P^3\times\mathbb P^1\hookrightarrow\mathbb P^7$,
and set $Y:=i_7(\mathbb P^3\times\mathbb P^1)$. By Example \ref{MinDeg} again, $\varepsilon(Y,\mathscr O_{\mathbb P^7}(1))=\frac{1}{2}$. The  computations are completely similar as in the previous case. Specifically, as it is easily seen, we have $c_1(\mathscr O_Y(1))=(P\times\mathbb P^1)+(\mathbb P^3\times y)$, $c_1(Y) =4(P\times\mathbb P^1)+2(\mathbb P^3\times y)$, $c_2(Y)= 6(\ell\times\mathbb P^1)+8(P\times y)$, 
$c_3(Y)$=$4(x\times\mathbb P^1)+12(\ell\times y)$, and $c_4(Y) = 8(x\times y)$, where $x\in\mathbb P^3$ and $y\in\mathbb P^1$ are points, $\ell$ is a line in $\mathbb P^3$ and $P$ is  a plane in $\mathbb P^3$. Using equation \eqref{chern3} to compute the Chern classes of $N$ we find:
\begin{eqnarray*}
 c_1(N)&=&4(P\times\mathbb P^1)+6(\mathbb P^3\times y), \\ c_2(N)&=& 6(\ell\times\mathbb P^1)+16(P\times y),\\
c_3(N)&=&4(x\times\mathbb P^1)+12(\ell\times y).\end{eqnarray*}
Using these formulae, \eqref{delta1}, \eqref{sc1} and \eqref{sc2} we get
$$\delta_{\eta}(Y,\mathscr O_{\mathbb P^7}(1))=20\eta(-1+9\eta-26\eta^2+30\eta^3-12\eta^4).$$
Setting $f(\eta)=-1+9\eta-26\eta^2+30\eta^3-12\eta^4$,  we have
$$\delta(Y,\mathscr O_{\mathbb P^7}(1))=\delta_{\frac{1}{2}}(Y,\mathscr O_{\mathbb P^7}(1))=10 f\big(\frac{1}{2}\big)=0,\;\;\text{but}\;\;\delta_{\frac{1}{3}}(Y,\mathscr O_{\mathbb P^7}(1))=\frac{20}{3} f\big(\frac{1}{3}\big)=\frac{20}{3}\times\frac{2}{27}>0.$$
Therefore $Y=i_7(\mathbb P^3\times\mathbb P^1)$ is $\mathscr O_{\mathbb P^7}(1)$-big.

\smallskip

Finally, consider at the Segre embedding  $i_8: \mathbb P^2\times\mathbb P^2\hookrightarrow\mathbb P^8$, and set $Y:=i_8(\mathbb P^2\times\mathbb P^2)$. In this case $Y=i(\mathbb P^2\times\mathbb P^2)$ is no more a subvariety of minimal degree because $\deg(Y)=6$ and $\codim_{\mathbb P^8}(Y)=4$. However, since $Y$ is not a linear subspace of 
$\mathbb P^8$ and $Y$ is still defined scheme-theoretically by quadratic equations, by Example \ref{MinDeg} we still have $\varepsilon(Y,\mathscr O_{\mathbb P^8}(1))=
\frac{1}{2}$. It is immediate to see that $c_1(\mathscr O_Y(1))=(\ell\times\mathbb P^2)+(\mathbb P^2\times\ell')$, $c_1(Y)=3(\ell\times\mathbb P^2)+3(\mathbb P^2\times\ell')
$, $c_2(Y)=3(x\times\mathbb P^2)+9(\ell\times\ell')+3(\mathbb P^2\times y)$, $c_3(Y)=9(\ell\times y)+9(x\times\ell')$, and $c_4(Y)=9(x\times y)$, where $x,y\in\mathbb P^2$ are points and  $\ell$ and $\ell'$ are lines in $\mathbb P^2$. Using again \eqref{chern3} and proceeding similarly as above we find for the Chern classes of $N$ the expressions:
\begin{eqnarray*}
c_1(N)&=&6(\ell\times\mathbb P^2)+6(\mathbb P^2\times\ell'),\\ c_2(N)&=&15(x\times\mathbb P^2)+27(\ell\times\ell')+15(\mathbb P^2\times y),\\
c_3(N)&=&45(\ell\times y)+45(x\times\ell'),\\ c_4(N)&=& 36(x\times y).\end{eqnarray*}
Now using \eqref{delta1}, \eqref{sc1} and \eqref{sc2} together with  the above relations,  we get (for $\eta\in(0,\frac{1}{2}$)):
$$\delta_{\eta}(Y,\mathscr O_{\mathbb P^8}(1))=9\eta^2(-10+80\eta-220\eta^2+252\eta^3-103\eta^4)=9\eta^2 f(\eta),$$
where $f(\eta):=-10+80\eta-220\eta^2+252\eta^3-103\eta^4$. Therefore $f(\frac{1}{2})=-10+40-55+31.5-6.4375=71.5-71.4375>0$, whence $\delta(Y,\mathscr O_{\mathbb P^8}(1))=\delta_{\frac{1}{2}}(Y,\mathscr O_{\mathbb P^8}(1))>0$. Thus $Y=i_8(\mathbb P^2\times\mathbb P^2)$ is $\mathscr O_{\mathbb P^8}(1)$-big.

\smallskip

Summing up we proved the following:

\begin{prop}\label{product} The images of the  Segre embeddings
$$i_5: \mathbb P^2\times\mathbb P^1\hookrightarrow\mathbb P^5, \;\;\; i_7:\mathbb P^3\times\mathbb P^1\hookrightarrow\mathbb P^7 \;\;\; and \;\;\; i_8: \mathbb P^2\times\mathbb P^2\hookrightarrow\mathbb P^8$$
are Seshadri $\mathscr O_{\mathbb P^k} (1)$-big, with $k=5$, $k=7$ and $k=8$ respectively.\end{prop}

\begin{rem*}\label{ofof} The examples of  Segre embeddings discussed above are  especially interesting. Indeed, although the rational normal scrolls $\mathbb F_e\hookrightarrow\mathbb P^{e+3}$ are all set-theoretic (but not scheme-theoretic) complete intersections in $\mathbb P^{e+3}$ (see \cite{Verdi}, or also \cite{BV}), the product $\mathbb P^2\times\mathbb P^1$ (respectively $\mathbb P^3\times\mathbb P^1$) is not set-theoretic complete intersection in $\mathbb P^5$ (respectively in $\mathbb P^7$). In fact, $\mathbb P^2\times\mathbb P^1$ (respectively  $\mathbb P^3\times\mathbb P^1$) is not even the zero locus of a section of an ample rank $2$ vector bundle on $\mathbb P^5$  (respectively of an ample rank $3$ vector bundle on $\mathbb P^7$), see
 \cite{Laz} for $\mathbb P^2\times\mathbb P^1$ and \cite[Corollary 4.5]{BV} for $\mathbb P^3\times\mathbb P^1$.
\end{rem*}
\end{prgrph*}
\begin{prgrph*}{\bf An elliptic surface}  The next  result provides an example of an irregular surface $Y$ in $\mathbb P^5$ which is
Seshadri $\mathscr O_{\mathbb P^5}(1)$-big.

\begin{prop}\label{elliptic} Consider the  geometrically ruled elliptic surface $Y=C\times\mathbb P^1$, with $C$ a smooth  elliptic plane curve, embedded in $\mathbb P^5$ via the composition  $C\times\mathbb P^1\hookrightarrow\mathbb P^2\times\mathbb P^1\hookrightarrow\mathbb P^5$, where the second inclusion is the Segre embedding of $\mathbb P^2\times\mathbb P^1$. Then $Y$ is Seshadri $\mathscr O_{\mathbb P^5}(1)$-big  in $\pn 5$.\end{prop}
\proof  Assume that the elliptic curve $C$ is  embedded in $\mathbb P^2$ as a cubic curve of equation $F(x_0,x_1,x_2)=0$. 
If $[y_0,y_1]\in \mathbb P^1$ then the homogeneous coordinates of a point of $\mathbb P^5$ are
$$z_{ij}=x_iy_j,\;\;\text{with $i=0,1,2$ and $j=0,1$}.$$
Thus $Y$ becomes an elliptic scroll in $\mathbb P^5$. Then by  \cite{SW} the homogeneous ideal of $Y$ in $\mathbb P^5$ is generated by the seven homogeneous polynomials $z_{10}z_{21}-z_{20}z_{11}$, $z_{20}z_{01}-z_{00}z_{21}$, $z_{00}z_{11}-z_{10}z_{01}$ and
$$\sum_{i=0}^2z_{is}\frac{\partial F}{\partial x_i}(z_{0t},z_{1t},z_{2t}),\;\;\text{where $0\leq s,t\leq 1$ . }$$
In particular, $Y$ is defined by $\mathscr O_{\mathbb P^5}(1)$ in degree $3$ in the terminology of \cite{BBF}. Then from \cite[Corollary 1.6]{BBF}  it follows
that $\varepsilon(Y,\mathscr O_{\mathbb P^5}(1))\geq\frac{1}{3}$. 

In our situation we have $\mathscr O_Y(1)=p_1^*(L)\otimes p_2^*(\mathscr O_{\mathbb P^1}(1))$, where $L=\mathscr O_C(1)$ is a line bundle of degree $3$ on $C$, and $\mathscr O_Y(K_Y)=p_2^*(\mathscr O_{\mathbb P^1}(-2)$. Using relation  \eqref{chern3} and computing as in the proof of Proposition \ref{Veronese2},  we find:
$$\deg(Y)=\mathscr O_Y(1)^{\cdot 2}=6,\;\; c_1(N)=c_1(p_1^*(L^{\otimes 6})\otimes p_2^*(\mathscr O_{\mathbb P^1}(4)), \;\;c_2(Y)=0,\;\;\text{and} \;\;\int_Yc_2(N)=54.$$  
These formulae yield
$$\int_Y(c_1(N)^{\cdot 2}-c_2(N))=90\;\;\;\text{and}\;\;\;\int_Yc_1(N)\cdot c_1(\mathscr O_Y(1))=30.$$ 
Therefore by \eqref{Ses2'} we find
$$\delta_{\eta}(Y,\mathscr O_{\mathbb P^5}(1))=\eta(-90\eta^2+90\eta-18).$$
Taking $\eta=\frac{1}{3}$ we finally get $\delta_{\frac{1}{3}}(Y,\mathscr O_{\mathbb P^5}(1))=\frac{-10+30-18}{3}=\frac{2}{3}>0.$  \qed
\end{prgrph*}

The following general result (which was noticed by Paoletti in \cite{Pao} when $(k,y)=(3,1)$) allows one to produce many new examples of Seshadri $A$-big or $A$-ample submanifolds, starting from some known ones.

\begin{prop}\label{finite} Let $f\colon X'\to X$ be a finite surjective morphism between smooth projective varieties. Let $A$ be an ample line bundle and set $A':=f^*(A)$. Let $Y$ be a smooth connected projective subvariety of $X$ and assume that $Y':=f^{-1}(Y)$ is smooth and connected. Then $\varepsilon(Y',A')=\varepsilon(Y,A)=:\varepsilon$ and $\delta_{\eta}(Y',A')=\deg(f)\,\delta_{\eta}(Y,A)$ for every $\eta\in(0,\varepsilon)$. In particular, $Y$ is Seshadri $A$-big $($respectively   $A$-ample$)$ if and only if $Y'$ is Seshadri $A'$-big  $($respectively $A'$-ample$)$.\end{prop}

\proof  Denote by $\mathscr I_Y$ the sheaf of ideals of $Y$ in $\mathscr O_X$. By a general result (see e.g. \cite[IV, (17.3.5)]{EGA}), as $f\colon X'\to X$ is a finite surjective morphism between  smooth varieties,  $f$ is a flat morphism. Therefore the canonical map $f^*(\mathscr I_Y)\to f^*(\mathscr O_X)=\mathscr O_{X'}$ is injective. It follows that $\mathscr I_{Y'}\cong f^*(\mathscr I_Y)$, where $\mathscr I_{Y'}$ is the sheaf of ideals of $Y'=f^{-1}(Y)$.    This shows that the blowing up $X'_{Y'}$ of $X'$ along $Y'$ identifies  with $\Proj(f^*(\oplus_{i=0}^{\infty}\mathscr I^i_Y))$ which, by \cite[II, Proposition 3.5.3]{EGA}, coincides with the fibered product $\Proj(\oplus_{i=0}^{\infty}\mathscr I^i_Y)\times_XX'=X_Y\times_XX'$. Then the finite surjective morphism $f$ induces by base-change a finite surjective morphism $\widetilde{f}\colon X'_{Y'}= X_Y\times_XX'\to X_Y$. As $\dim(X'_{Y'})=\dim(X_Y)$, $\widetilde f$ is also surjective. Clearly, $\widetilde f^*(E)=E'$, where $E'$ is the exceptional locus of $X'_{Y'}$. Then for every $\eta\in\mathbb R$,  the $\mathbb R$-line bundle $A+\eta E$ is ample if and only if the $\mathbb R$-line bundle $\widetilde f^*(A+\eta E)=A'+\eta E'$ is ample (because  $\widetilde{f}$ is finite and surjective). This fact and the definition of the Seshadri constant prove then the equality  $\varepsilon(Y',A')=\varepsilon(Y,A)=:\varepsilon$. 

Now by \eqref{delta1}, for every $\eta\in(0,\varepsilon)\cap\mathbb Q$ we have
\begin{equation}\label{delta112}\delta_{\eta}(Y,A)= -\sum_{t=0}^{k-2}{{k-2}\choose{t}}\eta^t\int_Ys_{y-k+t+2}(N)\cdot
c_1(A_Y)^{\cdot (k-t-2)}\end{equation} 
and
\begin{equation}\label{delta12}\delta_{\eta}(Y',A')= -\sum_{t=0}^{k-2}{{k-2}\choose{t}}\eta^t\int_{Y'}s_{y-k+t+2}(N')\cdot
{c_1(A'_{Y'})}^{\cdot (k-t-2)},\end{equation} 
where  $N$ and $N'$ are the normal bundles of $Y$ in $X$ and of $Y'$ in $X'$ respectively. If we set $g:=f|Y'\colon Y'\to Y$, then we have $N'=g^*(N)$, whence $s_{y-k+t+2}(N')=g^*(s_{y-k+t+2}(N))$ for every $t=0,1,\ldots,k-t-2$.
Since the morphism $g$ is flat (as a finite surjective morphism between smooth varieties), by \cite[Proposition 3.1, (d), p. 48]{Fu},  we get 
$$\int_{Y'}s_{y-k+t+2}(N')\cdot {c_1(A'_{Y'})}^{\cdot (k-t-2)}=\int_{Y'}g^*(s_{y-k+t+2}(N))\cdot c_1(g^*(A_Y))^{\cdot k-t-2}=$$
$$=\int_{Y'}g^*(s_{y-k+t+2}(N)\cdot c_1(A_Y)^{\cdot (k-t-2)})=\deg(g) \int_Ys_{y-k+t+2}(N)\cdot c_1(A_Y)^{\cdot (k-t-2)}, $$ 
for every $t=0,1,\ldots,k-2$. Substituting   in \eqref{delta112} and \eqref{delta12} and using the obvious equality $\deg(f)=\deg(g)$, we get $\delta_{\eta}(Y',A')=\deg(f)\,\delta_{\eta}(Y,A)$. This concludes the proof of the proposition.
\qed

\begin{rem*}\label{FulHan} If in Proposition \ref{finite} we take $X=\mathbb P^k$ then, by Fulton--Hansen connectedness theorem (see \cite{FH}), $f^{-1}(Y)$ is always connected.\end{rem*}

\begin{ex*}\label{Veroneseappl} To illustrate the use of Proposition \ref{finite}, let $Y$ be  the Veronese surface in $X=\mathbb P^5$. By Proposition \ref{Veronese2}, $Y$ is Seshadri $\mathscr O_{\mathbb P^5}(1)$-big. Fix an integer $n\geq 2$, and let $H_n$ be a smooth hypersurface of $\mathbb P^5$ that intersects $Y$ transversely (i.e.,  such that $C_n:=H_n\cap Y$ is a smooth curve). Consider the cyclic covering $f\colon X'\to X$ of degree $n$ branched along  $H_n\in|\mathscr O_{\mathbb P^5}(n)|$ and determined by $\mathscr O_{\mathbb P^5}(1)$ with $\mathscr O_{\mathbb P^5}(n)=\mathscr O_{\mathbb P^5}(H_n)$. As $H_n$ is smooth, $X'$ is also smooth. Moreover, by \cite[Lemma (17.1), p. 55]{BHPV}, the canonical class $\omega_{X'}$ is given by the formula 
$$\omega_{X'}=f^*(\omega_{\mathbb P^5}\otimes\mathscr O_{\mathbb P^5}(n-1))=f^*(\mathscr O_{\mathbb P^5}(n-7)).$$
It follows that  $X'$ is a Fano manifold for $n\leq 6$, $X'$ has trivial canonical class for $n=7$, and  the canonical class of $X'$ is ample for $n\geq 8$. Moreover, $Y'=f^{-1}(Y)$ is the cyclic covering $g\colon Y'\to\mathbb P^2$  of degree $n$  branched along the smooth curve $C_n\in|\mathscr O_Y(n)|=|\mathscr O_{\mathbb P^2}(2n)|$ (determined by the line bundle $\mathscr O_Y(1)=\mathscr O_{\mathbb P^2}(2)$ with $\mathscr O_Y(n)=\mathscr O_Y(C_n)$). As above, since $C_n$ is smooth, $Y'$ is also smooth. Moreover, again by \cite[Lemma (17.1), p. 55]{BHPV}, we have
$$\omega_{Y'}=g^*(\omega_{Y}\otimes\mathscr O_{Y}(n-1))=g^*(\mathscr O_{\mathbb P^2}(2n-5)).$$
In particular, if $n=2$,  $Y'$ is a Del Pezzo surface of degree $2$ (with canonical class $\omega_{Y'}=g^*(\mathscr O_{\mathbb P^2}(-1))$) while, if $n\geq 3$,  $Y'$ is a surface of general type (with ample canonical class $\omega_{Y'}=g^*(\mathscr O_{\mathbb P^2}(2n-5)$). As  $Y$ is Seshadri $\mathscr O_{\mathbb P^5}(1)$-big, $Y'$ is  Seshadri $f^*(\mathscr O_{\mathbb P^5}(1))$-big by  Proposition \ref{finite}.
\end{ex*}

\begin{rem*} Corollary \ref{Same0} of the next section will show that all the examples of Seshadri $\mathscr O_{\mathbb P^k}(1)$-big submanifolds of $\mathbb P^k$ given in this section are actually Seshadri $\mathscr O_{\mathbb P^k}(1)$-ample.
\end{rem*}

\section{Seshadri positivity and formal functions}\label{2a}
\addtocounter{subsection}{1}\setcounter{theorem}{0}

In this section we prove some general results which  give  further 
motivations for the study of Seshadri $A$-bigness and $A$-ampleness.
Let us start showing the following fact.

\begin{prop}\label{intersection} Under the above notation, if $Y$ is Seshadri $A$-ample in a smooth projective variety $X$
then $Y$ meets every irreducible hypersurface of $X$.\end{prop}

\proof Let $D$ be a closed irreducible  hypersurface of $X$, and
set $D'=\pi^{-1}(D)$. Then $D'$ is a hypersurface of $X_Y$. As  $\sh$ is the intersection of general very ample divisors $H_1,\ldots,H_{k-2} \in |mA^*-nE|$,  the intersection $D'\cap \sh$ is a curve in $\sh$. Finally, since $Y'$ is an ample divisor on the surface $\mathscr H$, the intersection  $Y'\cap D'=Y'\cap (D'\cap\sh)$  is not empty, which clearly  implies that $Y\cap D\neq\varnothing$. \qed

The proof of the first general result of the paper, Theorem \ref{Big2} below, is based upon the following lemma. To this purpose, recall that the {\em cohomological
dimension} of an algebraic scheme $V$ over the field $\mathbb C$ of complex numbers is 
$${\rm cd}(V):=\min\{t\;|\;H^i(V,\scrf)=0\; {\rm for \;any} \;i>t\;
{\rm and \;for\; any \;coherent \;sheaf} \;\scrf\; {\rm on}  \;V\}.$$

\begin{lemma}\label{Codlemma} Let $Y$ be a smooth $y$-dimensional subvariety  of a smooth
 $k$-dimensional  polarized variety $(X,A)$  such that $1\leq y\leq k-1$. If $Y$ is Seshadri $A$-ample then ${\rm cd}(X_Y\setminus Y')\leq
 k-2$.\end{lemma}
\proof  As usual, let  $H_1,\ldots,H_{k-2}$ be $k-2$ general members of $|mA^*-nE|$, 
 with $m$, $n$ big enough (so that $mA^*-nE$ is very ample). As above, set 
 $\sh:=H_1\cap\cdots\cap H_{k-2}$ and $Y':=\sh\cap E$.  For every $i=0,\ldots,k-3$, set
 $U_i:=(H_1\cap\cdots\cap H_{k-2-i})\setminus Y'$ and set $U_{k-2}:=X_Y\setminus
 Y'$. Then the closure
 $\overline{U}_i$ of $U_i$ is $\sh_i:=H_1\cap\cdots\cap H_{k-2-i}$ and
 $\overline{U}_{k-2}=X_Y$. It follows that $\overline{U}_i$ is a smooth
 connected projective
 variety of dimension $i+2$ for every $0\leq i\leq k-2$, because $H_1,\ldots,H_{k-2}\in|mA^*-nE|$
 are general.
 Moreover, every $\overline{U}_i$ contains $Y'$ and $U_{i-1}=U_i\cap
 H_{k-i-1}$, $1\leq i\leq
 k-2$. In other words $U_{i-1}$ is a hyperplane section of the (open)
 variety $U_i$. Then the conclusion  is a consequence of the following
 
 \begin{clam}\label{PortVale} ${\rm cd}(U_i)\leq i$ for every $0\leq i\leq
 k-2$.\end{clam}
 
 To prove the claim we proceed by induction on $i$. For $i=0$,
 $U_0=\sh\setminus Y'$, and since
 $Y$ is
 $A$-ample, $Y'$ is an ample effective divisor on $\sh$, whence $U_0$ is affine and
 hence ${\rm cd}(U_0)=0$
 by a well known affinity criterion of Serre.
 
 Assume therefore $1\leq i\leq k-2$ and set $U:=U_i$, $V:=U_{i-1}$ and
 $H:=H_{k-1-i}|U$. Then $\overline{U}=\sh_i$, $\overline{V}=\sh_{i-1}$,
 $H_{k-1-i}$ is a very ample divisor on  $\overline{U}$, and hence $H$ is a very
 ample divisor on $U$. By induction's hypothesis we may assume ${\rm cd}(V)\leq i-1$,
 and then we have to prove that ${\rm cd}(U)\leq i$. To do this, for every $p\in\mathbb Z$ consider the cohomology exact
 sequence
$$H^i(V,\so_V(pH))\to H^{i+1}(U,\so_U((p-1)H))\stackrel{\gra_p}{\to}H^{i+1}(U,\so_U(pH))
\to H^{i+1}(V,\so_V(pH))$$
 in which  the first and the last spaces are zero because
 ${\rm cd}(V)\leq i-1$. Therefore the maps $\gra_p$ are isomorphisms for all
 $p\in \zed$.
 Then by \cite[Proposition 3.1, p. 92]{Ample}, the inequality  ${\rm cd}(U)\leq i$ is
 equivalent to the vanishings $H^{i+1}(U,\so_U(pH))=0$ for $p\ll 0$, or else, by using the
 isomorphisms $\gra_p$, to the vanishings $H^{i+1}(U,\so_U(pH))=0$ for $p\gg 0$. To prove these
 latter vanishings consider the exact sequence of local cohomology
 $$H^{i+1}(\overline{U},\so_{\overline{U}}(pH_{k-1-i}))\to H^{i+1}(U,\so_U(pH))
 \to H^{i+2}_{Y'}(\overline{U},\so_{\overline{U}}(pH_{k-1-i})).$$ The first
 space is zero for $p\gg 0$ by Serre's theorem, so that it will be enough to
 show that the last space is also zero for $p\gg 0$. By using formal duality \cite[Theorem 3.3, p. 94]{Ample} and the fact that $\dim(\overline{U})=i+2$, we get
 $$H^{i+2}_{Y'}(\overline{U},\so_{\overline{U}}(pH'))^*\cong
 H^0({\overline{U}}_{/Y'},
 [{\omega_{\overline{U}}(-pH')}]_{/Y'}),$$ 
 where $H':=H_{k-1-i}$,
 $\omega_{\overline{U}}(-pH')=\omega_{\overline{U}}\otimes\so_{\overline{U}}(-pH')$, and $\scrf_{/Y'}$ denotes the formal completion along $Y'$ of a sheaf $\scrf$ in $\overline{U}$. So everything
 reduces to show that the latter formal cohomology space vanishes for $p\gg 0$.
 
 First we observe that the normal bundle $N_{Y'|\overline{U}}$ of $Y'$ in
 $\overline{U}$ is of the form
 \begin{equation}\label{decomp}
 N_{Y'|\overline{U}}\cong M\oplus L^{\oplus i},
 \end{equation}
 where $M=N_{Y'|\sh}$ and $L=\so_{Y'}(mA^*-nE)$. As the divisor $mA^*-nE$ is very ample, the line bundle $L$ is very ample. Moreover, $M$ is ample because by hypothesis $Y$ is $A$-ample. To prove \eqref{decomp}, it is enough to observe that in $\overline{U}=H_1\cap\cdots\cap H_{k-2-i}$, the curve $Y'$ is the
 proper intersection
 of $\sh$ with $\overline{E}:=E\cap\overline{U}$, whence
 $$N_{Y'|\overline{U}}\cong N_{Y'|\sh}\oplus N_{Y'|\overline{E}}=
  M\oplus N_{\sh|\overline{U}}|Y'\cong M\oplus L^{\oplus i}.$$
 
 Now, for every $\ell\geq 0$, let
 $Y'_\ell=(Y',\so_{\overline{U}}/\sy_{Y'}^{\ell +1})$ be the
 $\ell$-th infinitesimal neighbourhood of $Y'$ in $\overline{U}$, where
 $\sy_{Y'}$ is the ideal
 sheaf of $Y'$ in $\overline{U}$. Then one has
 \begin{equation}\label{formal1}H^0({\overline{U}}_{/Y'},
[{\omega_{\overline{U}}(-pH')}]_{/Y'})=\varprojlim_{\ell} H^0(Y'_{\ell},
 \omega_{\overline{U}}\otimes\so_{\overline{U}}(-pH'){|Y'_{\ell}}).\end{equation}
 Moreover the standard exact sequence 
$$0\to\sy_{Y'}^{\ell +1}/\sy_{Y'}^{\ell +2}\cong{\bf S}^{\ell
 +1}(N_{Y'|\overline{U}}^{\vee})\to \so_{Y'_{\ell +1}}\to\so_{Y'_\ell}\to 0$$
 shows that for every $\ell\geq 0$ the kernel of the map
 $$\varepsilon_\ell:H^0(Y'_{\ell+1},\omega_{\overline{U}}(-pH'){|Y'_{\ell+1}})
 \to H^0(Y'_{\ell},\omega_{\overline{U}}(-pH')|Y'_{\ell})$$ 
 is $H^0(Y',{\bf S}^{\ell+1}(N_{Y'|\overline{U}}^{\vee})\otimes\omega_{\overline{U}}(-pH')|Y')$. We claim that
 \begin{equation}\label{PortVale1} H^0(Y',{\bf S}^{\ell+1}(N_{Y'|\overline{U}}^{\vee})\otimes\omega_{\overline{U}}(-pH')|Y')=0,\;\text{for all $p\gg 0$ and all $\ell\geq 0$}.\end{equation}
 Indeed, choose an integer $p_0>0$ such that $\omega_{\overline{U}}(-pH')|Y'$ is ample for every $p\geq p_0$. Since $ {\bf S}^{\ell+1}(N_{Y'|\overline{U}}^{\vee})\otimes\omega_{\overline{U}}(-pH')|Y'$ is a direct sum of bundles of the form $M^{-\alpha}\otimes L^{-\beta}\otimes\omega_{\overline{U}}(-pH')|Y'$ (whose inverses are ample) we have
$$H^0(Y',M^{-\alpha}\otimes L^{-\beta}\otimes\omega_{\overline{U}}(-pH')|Y')=0,\;\text{for every $p\geq p_0$,}$$
which proves the claimed assertion (\ref{PortVale1}).  It  thus follows that the maps $\varepsilon_\ell$ of the projective system are injective for every $\ell\geq 0$ and $p\geq p_0$. Then from \eqref{formal1} it follows that 
 $$H^0({\overline{U}}_{/Y'},[{\omega_{\overline{U}}(-pH')}]_{/Y'})\cong H^0(Y',\omega_{\overline{U}}(-pH'){|Y'}),\;\; \text{for all $p\geq p_0$}.$$
 In particular, the vanishing of $H^0({\overline{U}}_{/Y'},{[\omega_{\overline{U}}(-pH')]}_{/Y'})$ for $p\geq p_0$ is equivalent to the vanishing of  $H^0(Y',\omega_{\overline{U}}(-pH'){|Y'})$ for $p\geq p_0$. Finally, since $\dim(Y')=1>0$, the latter vanishing is  obvious, so that Claim \ref{PortVale} is proved. This completes the proof of Lemma \ref{Codlemma}. \qed

Theorem \ref{Big2} below will show that the Seshadri $A$-bigness and the Seshadri $A$-ampleness of
$Y\subset X$ are related to properties of formal rational functions along
$Y$ and to the cohomological dimension of the complement $U:=X\setminus Y$.
First recall some definitions we need. We refer to \cite[Chapters III,
V]{Ample}, or also to \cite{B2}, for more details.

\begin{deff*}\label{FF}Let $V$ be an integral algebraic scheme over $\mathbb C$, $Y\subset V$ a closed subscheme,
$V_{/Y}$ the formal completion of $V$ along $Y$, and let $K(V_{/Y})$
be the ring of the formal rational functions of $V$ along $Y$ (see \cite{HM},
or also \cite{Ample}). Then $K(V_{/Y})$ is a field if $V$ is
non-singular and $Y$ is connected (loc. cit.). Thus in this case there is a natural homomorphism of $\mathbb C$-algebras
$K(V)\to K(V_{/Y})$ making $K(V_{/Y})$ a field extension of $K(V)$. According to \cite{HM} we say that
\begin{enumerate}
\item[{\rm i)}] $Y$ is ${\rm G}2$  in $V$ if the map $K(V)\to K(V_{/Y})$
makes $K(V_{/Y})$ a finite field extension of $K(V)$;
\item[{\rm ii)}] $Y$ is ${\rm G}3$  in $V$ if the map $K(V)\to
K(V_{/Y})$ is an isomorphism.\end{enumerate}
\end{deff*}

\begin{theorem}\label{Big2} Let $Y$ be a smooth subvariety of dimension $y$ of a smooth
complex projective polarized variety $(X,A)$ of dimension $k$ such that $1\leq y\leq k-1$. Then the following statements hold:
\begin{enumerate}
\em\item\em If $Y$ is Seshadri $A$-big then $Y$ is ${\rm G}2$ in $X$ and $k-y-1\leq\cd(X\setminus Y)\leq k-1$.
\em\item\em If $Y$ is Seshadri $A$-ample then $Y$ is ${\rm G}3$ in $X$ and $k-y-1\leq\cd(X\setminus Y)\leq k-2$.\end{enumerate}\end{theorem}

\proof Under the usual notation, let $Y':=E\cap\mathscr H$,
where  $\sh:=H_1\cap\cdots\cap H_{k-2}$ is the intersection of
$k-2$ general very ample divisors $H_1,\ldots,H_{k-2}\in|mA^*-nE|$.
As the intersection $Y'=E\cap \mathscr H$  is proper and
$N_{\sh|X_Y}|Y'\cong N_{Y'|E}$, the standard exact sequence of normal
bundles  
$$0\to N_{Y'|\mathscr H}\to  N_{Y'|X_Y} \to N_{\mathscr H|X_Y}|Y'\to 0$$
 splits to give
$$ N_{Y'|X_Y}\cong N_{Y'|\sh}\oplus N_{Y'|E}.$$
The normal bundle $N_{Y'|\sh}$ is ample by definition of Seshadri $A$-bigness. On
the other hand $N_{Y'|E}$ is also ample since $Y'$ is a complete intersection
in $E$. It follows that $N_{Y'|X_Y}$ is ample. From this and a result of
Hartshorne (see \cite[p. 198]{Ample}, or also \cite[\S 6]{CDAV}) it follows
that $Y'$ is ${\rm G}2$ in $X_Y$.

Now, as the morphism $\pi$ is proper and birational and $E=\pi ^{-1}(Y)$, by a theorem of Hironaka and Matsumura \cite{HM}
we get
\begin{equation}\label{formal}
K({X}_{/Y})\cong K({X_Y}_{/E}),
\end{equation}
where ${X}_{/Y}$ (respectively ${X_Y}_{/E}$) denotes the formal
completion of $X$ (respectively of $X_Y$) along $Y$ (respectively along $E$).
Since $Y'\subset E$, we get natural maps
$$K(X_Y)\to K({X_Y}_{/E})\to K({X_Y}_{/Y'}),$$
where ${X_Y}_{/Y'}$ denotes the formal completion of $X_Y$ along
$Y'$. The fact that $Y'$ is ${\rm G}2$ in $X_Y$ means that
$K({X_Y}_{/Y'})$ is a finite field extension of $K(X)=K(X_Y)$ (and in
particular, $K({X_Y}_{/Y'})$ is a field). On the other hand,
$K({X_Y}_{/E})$ is a field because $X_Y$ is smooth and $E$ is
connected. In particular, the second map is injective, and so
$K({X_Y}_{/E})\cong K({X}_{/Y})$ is a subfield of
$K({X_Y}_{/Y'})$, whence a finite field extension of $K(X_Y)\cong
K(X)$. In other words, $Y$ is ${\rm G}2$ in $X$.

 To prove $1)$ observe that, since $\dim(Y)=y$  and $\dim(X)=k$, a general complete intersection of $X$ of
 codimension $k-y-1$ does not meet $Y$ and hence ${\rm cd}(X\setminus Y)\geq
 k-y-1$. Moreover,
 as $Y\neq \varnothing$, a theorem of Hartshorne and Lichtenbaum (see
 \cite{CDAV} or \cite[p. 98]{Ample}) implies that ${\rm cd}(X\setminus
 Y)\leq k-1$. So the first statement   is proved. 

To prove $2)$,  assume   that $Y$ is $A$-ample. Then, by Lemma \ref{Codlemma}, $\cd(X_Y\setminus Y')\leq
 k-2$ and, since $Y'$ is G$2$ in $X_Y$,  a result of Speiser  applies to say  that $Y'$ is in fact G$3$ in $X_Y$, i.e., the field extension $K(X_Y)\subseteq K({X_Y}_{/Y'})$ is an isomorphism (see \cite{Sp}, or also \cite[Corollary 2.2, p. 202]{Ample}).
As $Y'\subset E$, $K({X_Y}_{/E})$ is a subfield extension of $K({X_Y}_{/Y'})$, whence the field extension $K(X)=K(X_Y)\subseteq K({X_Y}_{/E})=K(X_{/Y})$
is also an isomorphism, i.e., $Y$ is ${\rm G}3$ in $X$. Since  $Y$ meets every hypersurface of $X$  by Proposition \ref{intersection} and since $Y$ is ${\rm G}3$ in $X$, we finally get $\cd(X\setminus Y)\leq k-2$ by Speiser's result again. Thus part $2)$ also holds, and this completes the proof. \qed

 \section{Comparing Seshadri bigness and Seshadri ampleness}\label{2b}
\addtocounter{subsection}{1}\setcounter{theorem}{0}

Now we come up to the second general main result of this paper (which generalizes Theorem 3.1 of \cite{BBF} to the case when $Y$ is a submanifold of dimension $\geq 2$).  Note that Theorem \ref{Main} below is new even in the case when $\codim_X(Y)=1$. We keep the notation as in Section \ref{def}.

\begin{theorem}\label{Main} Let $(X,A)$ be a polarized manifold of
dimension $k\geq 3$, and let $Y$ be a
 submanifold in $X$ of dimension $y\geq 1$ which is Seshadri $A$-big.
Then either $Y$ is $A$-ample, or there exists an irreducible hypersurface $D$ of $X$ such that
 $Y\cap D=\varnothing$. In the latter case the set of all irreducible hypersurfaces $D$ of $X$ such that $Y\cap D=\varnothing$ is
 finite.\end{theorem}
 
 \proof By the $A$-bigness assumption on $Y$,  we can fix an $\eta=\frac{n}{m}\in(0,\gre(Y,A))\cap \rat$ such that
 $\delta_{\eta}(Y,A)>0$, and  we may assume
  the linear system $|mA^*-nE|$ to be  very ample. Let $i:X_Y\hookrightarrow P:={\Bbb
 P}(H^0(\so_{X_Y}(mA^*-nE))^{\vee})$ be the corresponding closed embedding, and let
 $L$ be a general $(k-1)$-codimensional linear subspace of $P$. Then the linear system
 of all $(k-2)$-codimensional linear subspaces of $P$ containing $L$ is
 parameterized by the $(k-2)$-dimensional projective space $\pn {k-2}$. For every
 $t\in \pn {k-2}$, denote by $\sh'_t$ the $(k-2)$-codimensional linear subspace
 of $P$ containing $L$ and corresponding to $t$, and let  $\Lambda:=\{\sh'_t\}_{t\in\pn {k-2}}$. Setting
 $\sh_t:=\sh'_t\cap X_Y$, consider the incidence correspondence of  the family $\Lambda$,
 $$V':=\{(t,x)\in\pn {k-2}\times X_Y\;|\;x\in\sh_t\},$$ and denote by $p:V'\to\pn
 {k-2}$ the restriction of the first projection of $\pn {k-2}\times X_Y$, and by
 $\vphi:V'\to X_Y$ the restriction of the second projection. Note that, by
 construction, $\vphi$ coincides with the blowing up of $X_Y$ along $L\cap
 X_Y=\sh_\gra\cap\sh_\grb$, where $\gra$ and $\grb$ are any two different
 points of
 $\pn {k-2}$. In particular $V'$ is smooth because $X_Y$ is smooth and the
 linear subspace $L$ is general (and so,  by Bertini's theorem, $L\cap X_Y$ is
 smooth). Let $W'$ be the proper inverse image of $E$ under $\vphi$, i.e., the
 blowing up of $E$ along $L\cap E$. Since the linear subspace $L$ is general, again by Bertini's theorem,
$E\cap L$ is a finite set of reduced points. In particular, $W'$ is also smooth because $E$ is so.
Then we get the commutative diagram
 $$
  \xymatrix{ W' \  \ar@{^(->}[rr]^{} \ar[rd]_{q} &  & V' \ar[ld]^{p} \\ & \pn {k-2} & } 
  $$
 with $p$ and $q$ proper morphisms. Note that for every $t\in \pn {k-2}$,
 $p^{-1}(t)\cong\sh_t$
 and $q^{-1}(t)\cong E\cap\sh_t$ under $\vphi$. Set $Y_t:=q^{-1}(t)=E\cap \sh_t$.
 
 We claim that the morphisms $p$ and $q$ are both flat. Indeed they are
 proper surjective
 morphisms between smooth varieties and their fibers have constant dimension
 ($2$ and $1$, respectively). Then the assertion follows from a well known criterion of flatness
 due to Grothendieck,
 see \cite{SGA}, or also \cite[V, (3.5)]{AK}.
 
 As $Y$ is Seshadri $A$-big and $L$ is general, the normal bundle of $q^{-1}(t)=Y_t$ in
 $p^{-1}(t)=\sh_t$ is ample for $t\in\pn {k-2}$ general. Set
 $$B:=\{t\in \pn {k-2} \;|\; N_{Y_t|\sh_t}\; {\rm is \;ample}\},\;\;
 V:=p^{-1}(B)\;\;{\rm  and}\; W:=q^{-1}(B).$$
Since the ampleness of  $N_{Y_t|\sh_t}$ is an open condition (see \cite{EGA}, III), $B$ is an open non-void  subset of $\pn {k-2}$ and we get the commutative
 diagram
$$ \xymatrix{ W \  \ar@{^(->}[rr]^{} \ar[rd]_{q} &  & V \ar[ld]^{p} \\ & B & } $$
(where, by a slight abuse of notation, the restrictions $q|W$ and $p|V$ are still denoted by $q$ and $p$ respectively).
 Note that, since $p$ is flat, $N_{W|V}|Y_t\cong N_{Y_t|\sh_t}$ for every
 $t\in B$. Therefore
 $N_{W|V}$ is a $q$-ample line bundle where, as usual, $N_{W|V}$ denotes the
 normal bundle of $W$ in $V$.
 \begin{clam}\label{spanned} There exists a commutative diagram
 $$
 \xymatrix{ V \  \ar[rr]^{f} \ar[rd]_{p} &  & V^* \ar[ld]^{p^*} \\ & B & } 
 $$
 where $V^*$ is a normal variety, $p^*$ is a proper surjective morphism, and
 $f$ is a proper birational morphism which is an isomorphism in a
 neighbourhood of $W$ {\rm (}and, in particular, $W$ is embedded in $V^*$ as
 a Cartier divisor{\rm )} such that  $\so_{V^*}(W)$ is $p^*$-ample. \end{clam}
 
 The claim is proved in \cite[III, Theorem 4.2]{Ample} in the case when $B$
 is a point, and
 extended to the general case in \cite{Ba}. Note that in \cite{Ba} the
 assumption that $V$ and $W$
 were projective played no role; the only thing used in the proof was that
 the morphisms $p$ and
 $q$ were proper. Note also (see \cite{Ba}) that the morphism $f$ is
 gotten in the following
 canonical way. One checks  first that for $r\gg 0$ the natural map
 $p^*p_*(\so_V(rW))\to\so_V(rW)$ is surjective. Therefore for $r\gg 0$ there is a unique
 $B$-morphism $g:V\to {\Bbb P}(p_*(\so_V(rW)))$ and one shows that one can take as $f$ the morphism with
 connected fibers
 arising from the Stein factorization of $g$.
 
 Since $V^*$ is normal the general fiber of $p^*$ is also normal. Therefore,
 shrinking $B$ if
 necessary,  we may assume that all the fibers of $p^*$ are normal. It
 follows that for every
 $t\in B$ the restriction $f_t:=f |p^{-1}(t): \sh_t\to {p^*}^{-1}(t)$
 is obtained by taking the Stein factorization of the morphism defined by
 $|\so_{X_Y}(rE)\otimes\so_{\sh_t}|$ for $r\gg 0$,  so in particular $f_t$
 depends only on $\sh_t$, and not on the choice of the linear subspace $L$ of
 codimension $k-1$ contained in $\sh'_t$.
 
 Let $Z$ (respectively $Z_t$) be the locus of $V$ (respectively of
 $p^{-1}(t)=\sh_t$) consisting
 of all points at which the morphism $f$ (respectively $f_t$) is not a
 biregular isomorphism. Then
 it is easy to see that $Z\cap p^{-1}(t)=Z_t$ for every $t\in B$.
 
 Now assume  that $Y$ is not $A$-ample. Then using the above observations it
 follows easily that
 there exists an irreducible component $Z^*$ of $Z$ such that $\dim(Z^*\cap
 p^{-1}(t))\geq 1$
 for every
 $t\in B$, whence
 $$\dim(Z^*\cap p^{-1}(t))=1\;{\rm for\; every}\; t\in B,$$
 because $p$ is of relative dimension $2$. As $\dim(B)=k-2$, we have $\dim(Z^*)=k-1$.
 Let $Z'$ be the closure of $Z^*$ in $V'$. We claim that
 \begin{equation}\label{pidim}
 \dim(\pi(\vphi(Z'))=k-1.\end{equation}
 To prove (\ref{pidim}), recall that by Claim \ref{spanned} the morphism
 $f$ is an isomorphism
 in a neighbourhood of $W$, whence $W\cap Z^*=\varnothing$. It thus follows that
 $Y_t\cap Z'=Y_t\cap
 Z^*=\varnothing$ for all $t\in B$. Therefore $Z'$ could intersect $W'$ at
 most at points belonging
 to  the union of all fibers $q^{-1}(t)=Y_t$ with $t\in \pn {k-2}\setminus B$. In
 particular, $Z'$ is
 not contained in $W'$, and hence $\vphi(Z')$ is not contained in $E$.
 Therefore, to prove
 (\ref{pidim}), it will be enough to show that $\dim(\vphi(Z'))=k-1$
 (because $\pi$ is an  isomorphism off $E$).
 
 Assume that this last equality does not hold, i.e., $\dim(\vphi(Z'))< k-1$.
 Recalling that
 $\vphi$ is the blowing up of $X_Y$ along $L\cap X_Y= \sh_\gra\cap\sh_\grb$,
 with two
 fixed different points $\gra$, $\grb\in B$, it follows that
 $\vphi(Z')=L\cap X_Y$. As noted
 above, the morphism $f_\gra=f{|\sh_\gra}$ does not change if we replace the
 linear subspace $L$ of $P$ by another $(k-1)$-codimensional  linear subspace $M$
 of $P$ ($L\neq M$) of the form $M=\sh_\gra\cap\sh''$, where $\sh''$ is a
 $(k-2)$-codimensional linear subspace of $P$, which  does not belong
 to the family $\Lambda=\{\sh'_t\}_{t\in\pn {k-2}}$. In other words,
 $f{|\sh_\gra}$ does not change if we replace the family
 $\Lambda$ by a different family
 $\Omega=\{\sk'_s\}_{s\in \pn {k-2}}$ of $(k-1)$-codimensional linear subspaces
 of $P$, such that $\sh'_\gra=\sk'_\gamma$ for some $\gamma\in \pn {k-2}$. Thus
 we can assume that there is an irreducible curve $C$ of $\sh_\gra$ contracted
 to a point by $f_\gra=f{|\sh_\gra}$ and which is not contained in $M\cap X_Y$.
 In other words, varying $L$ a bit, we can assume that $\vphi(Z')$ is not
 contained in $L\cap X_Y$, which proves equality (\ref{pidim}).
 
 To prove the first part of the theorem, it is sufficient to show
 that \begin{equation}\label{empty}
 \vphi(Z')\cap E=\varnothing.
 \end{equation}
 Indeed, from (\ref{pidim}) and (\ref{empty}) it follows that $D:=\pi(\vphi(Z'))$
 is a hypersurface of $X$ such that $D\cap Y=\varnothing$.
 
  Assuming that (\ref{empty}) fails
 to be true, pick a point $x\in\vphi(Z')\cap E$. Then, since $X_Y$ is smooth and recalling equality (\ref{pidim}) , we
 have
 $$\dim(\vphi(Z')\cap E)\geq
 \dim(\vphi(Z'))+\dim(E)-\dim(X_Y)=\dim(\vphi(Z'))-1=k-2.$$
 Therefore there is an irreducible subvariety $\Gamma\subseteq\vphi^{-1}(Z')\cap
 E$ of dimension $k-2$ passing through $x$. As $\sh_t$ is a
 $(k-2)$-dimensional complete intersection in $X_Y$, 
 $\sh_t\cap\Gamma\neq \varnothing$ for every $t\in \pn {k-2}$. In particular, for
 every $t\in B$,
 $$(Z^*\cap W)\cap p^{-1}(t)=(Z'\cap W')\cap p^{-1}(t)\neq \varnothing,$$
 which contradicts the fact that $Z^*\cap W=\varnothing$. This proves
 (\ref{empty}), and thereby the first part of the theorem.
 
 To prove the second part of the statement, let $D$ be a hypersurface of $X$
 such that $Y\cap D=\varnothing$. Set $D':=\pi^{-1}(D)$. Then $D'\cap
 E=\varnothing$. Let $D^*$ be the proper transform of $D'$ via $\vphi:V'\to X_Y$.
 If $D^*\cap V=\varnothing$ then $p(D^*)\subseteq \pn {k-2}\setminus B$, whence
 $\dim(p(D^*))\leq k-3$. On the other hand since $\dim(D^*)=k-1$ and $p$ is of
 relative dimension $2$, one gets $\dim(p(D^*))\geq k-3$, whence
 $\dim(p(D^*))=k-3$. It follows that $D^*$ is an irreducible component of the
 algebraic subset $p^{-1}(\pn {k-2}\setminus B)$ of $V$, so that there are
 finitely many
 possibilities for $D^*$ (and hence also for $D$) if $D^*\cap V=\varnothing$.
 Assume now that $D^*\cap V\neq\varnothing$. As $D'\cap E=\varnothing$,  we get
 $D^*\cap W'=\varnothing$, whence $\widetilde{D}:=D^*\cap V$ does not meet $W$.
 Since $\dim(\widetilde{D})=k-1$, $\widetilde{D}\cap p^{-1}(t)$ is a curve on
 $p^{-1}(t)=\sh_t$ which does not meet $Y_t$ for all $t\in B$. Since $f_t$ is
 given on $p^{-1}(t)$ by the linear system $|rY_t|$, $r\gg 0$, it follows that
 $\widetilde{D}\cap p^{-1}(t)\subseteq Z_t$ for every $t\in B$, whence
 $\widetilde{D}\subseteq Z$. If fact, for dimension reasons,
 $\widetilde{D}$ is an
 irreducible component of $Z$, and there are only finitely many irreducible
 components of $Z$. This completes the proof of the theorem. \qed
 
\begin{corollary}\label{Same0} Let $X$ be a $k$-dimensional projective manifold,
$k\geq 3$, polarized by an ample line bundle $A$. Let $Y$ be a smooth
subvariety of $X$ of dimension $y\geq 1$. Then the following statements hold:
\begin{enumerate}
\em\item\em $Y$ is Seshadri $A$-ample if and only if $Y$
is Seshadri $A$-big and $\cd(X\setminus Y)\leq k-2$.
\em\item\em $Y$ is Seshadri  $A$-ample if and only if $Y$
is Seshadri $A$-big
and $Y$ meets every hypersurface $D$ in $X$.\end{enumerate}

In particular, for $X=\pn k$,  Seshadri $\mathscr O_{\mathbb P^k}(1)$-bigness  is equivalent to is Seshadri $\mathscr O_{\mathbb P^k}(1)$-ampleness.
\end{corollary}

\proof It follows immediately from Proposition \ref{intersection}, Theorem \ref{Big2} and Theorem
\ref{Main}. \qed

\begin{corollary}\label{div2} Let $Y$ be a submanifold of codimension one of a polarized manifold $(X,A)$ which is Seshadri $A$-big $($e.g. if the normal bundle $N$ of $Y$ in $X$ satisfies the condition \eqref{22} of Corollary $\ref{codim111})$. Then the set of irreducible hypersurfaces of $X$ that do not meet $Y$ is finite.\end{corollary}

\proof Everything follows from Corollary \ref{codim111} and Theorem \ref{Main}. \qed

\begin{ex*}\label{exceptional} Let $Y$ be a smooth subvariety of $X=\mathbb P^k$ of dimension $y\geq 1$. By Corollary \ref {Same0}, $Y$ is $\mathscr O_{\mathbb P^k}(1)$-ample if and only if it is Seshadri $\mathscr O_{\mathbb P^k}(1)$-big. Now, assume $Y$ is $\mathscr O_{\mathbb P^k}(1)$-ample in $\mathbb P^k$, and let $x\in\mathbb P^k\setminus Y$. Let $u:X'\to \pn k$ be the blowing up of $\mathbb P^k$ at $x$, and denote by  $F=u^{-1}(x)$ the exceptional locus of $X'$. Then the closed embedding $i\colon Y\hookrightarrow \mathbb P^k$ lifts to a (unique) closed embedding $Y\hookrightarrow X'$. As the line bundle $\mathscr O_{X'}(-F)$  is $u$-ample, there is a sufficiently big $n\in\mathbb N$ such that $A:=u^*(\mathscr O_{\mathbb P^k}(n))\otimes \mathscr O_{X'}(-F)$ is ample on $X'$. Since $Y\cap F=\varnothing$, $A_Y\cong u^*(\mathscr O_{\mathbb P^k}(n))|Y\cong \mathscr O_{Y}(n)$. Then we claim that $Y$ is Seshadri $A$-big, but not $A$-ample. Indeed, the fact that $Y$ is not $A$-ample  follows from Corollary \ref{Same0}$(2)$,  because $Y\cap F=\varnothing$. On the other hand, since $Y$ is Seshadri $\mathscr O_{\mathbb P^k}(1)$-big in $\mathbb P^k$, by Remark \ref{indep} (and especially by \eqref{indep1}), $Y$ is also Seshadri $\mathscr O_{\mathbb P^k}(n)$-big in $\mathbb P^k$, i.e.,
\begin{equation}\label{delta6}\delta_{\eta}(Y,\mathscr O_{\mathbb P^k}(n))=-\sum_{t=0}^{k-2}{{k-2}\choose{t}}\eta^t\int_Ys_{y-k+t+2}(N_{Y|\mathbb P^k})\cdot
c_1(\mathscr O_Y(n))^{\cdot (k-t-2)}>0, \end{equation}
for some $\eta\in (0,\varepsilon(Y,\mathscr O_{\mathbb P^k}(n))\cap\mathbb Q$.
Since $Y\cap F=\varnothing$ and  $N_{Y|X'}=N_{Y|\mathbb P^k}$, we also have, by  \eqref{delta6},
\begin{eqnarray*}
\delta_{\eta}(Y,A)&=&-\sum_{t=0}^{k-2}{{k-2}\choose{t}}\eta^t\int_Ys_{y-k+t+2}(N_{Y|X'})\cdot c_1(A_Y)^{\cdot (k-t-2)}\\&=&-\sum_{t=0}^{k-2}{{k-2}\choose{t}}\eta^t\int_Ys_{y-k+t+2}(N_{Y|\mathbb P^k})\cdot
c_1(\mathscr O_Y(n))^{\cdot (k-t-2)}\\&=&\delta_{\eta}(Y,\mathscr O_{\mathbb P^k}(n))>0.\end{eqnarray*}

Moreover, we could blow up a smooth subvariety of $\mathbb P^k$ which does not intersect $Y$ (instead  a point $x\in\mathbb P^k\setminus Y$), and the same conclusion as above would be  true.
\end{ex*}

Example \ref{exceptional} suggests that the (finitely many) irreducible hypersurfaces of $X$ that do not intersect a given Seshadri $A$-big submanifold $Y$ of $X$ should be rather ``special''. The next result gives some evidence in this sense.

 \begin{theorem}\label{Intersection} Let $(X,A)$ be a smooth projective variety of
 dimension $k\geq 3$ polarized by an ample line bundle $A$, and let $Y$ be a
smooth closed subvariety of $X$ of dimension $y\geq 1$ which is Seshadri $A$-big. Let $Z$ be an irreducible
 hypersurface of $X$ which is a local complete intersection in $X$.
 If $Y$ is Seshadri $A$-big and the normal bundle $N_{Z|X}$
 of $Z$ in $X$ is ample, then $Y\cap Z\neq \varnothing$.\end{theorem}
 
 \proof  Because the normal bundle $N_{Z|X}$ is ample, by  \cite[Theorem 4.2, p. 110]{Ample} there exists a birational morphism $\varphi\colon X\to X_1$ onto a normal projective variety $X_1$, which is a biregular isomorphism between an open neighbourhood of  $Z$ in $X$ and  an open neighbourhood of $\varphi(Y)$ in $X_1$, such that $\varphi(Z)$ becomes an ample effective divisor on $X_1$.
 As $Y$ is Seshadri $A$-big, by Theorem \ref{Big2}$(1)$, $Y$ is G$2$ in $X$. Then by an observation of Hironaka and Matsumura (\cite[top of p. 65]{HM}), $\dim\varphi(Y)>0$. Finally, as $\varphi(Z)$ is an effective ample divisor on $X_1$, $\varphi(Y)\cap\varphi(Z)\neq
\varnothing$, which, by Zariski's Connectedness Theorem, see \cite[III, Th\'eor\`eme (4.3.1), p. 130]{EGA},  implies that $Y\cap Z\neq\varnothing$. \qed

The last result shows that Seshadri $A$-ampleness and Seshadri $A$-bigness are open properties under smooth deformations. We omit the proof since it is almost identical to that of  \cite[Theorem 4.1]{BBF}. 

 \begin{theorem}\label{Deform} Let $X$ be a smooth projective variety and
 let $f:X\to B$ be a
 smooth proper morphism of relative dimension $k\geq 3$ onto an irreducible
 algebraic variety $B$
 of positive dimension. Let $Y$ be a closed subvariety of $X$ such that the
 restriction
 $g:=f|Y\colon Y\to B$ is a smooth morphism of relative dimension $y\geq 1$. Let
$A$ be an $f$-ample line bundle on $X$. For any $t\in B$ set $X_t:=f^{-1}(t)$, $Y_t:=g^{-1}(t)$ and
 $A_t:=A|X_t$. Assume that there exists a point $t_0\in B$ such that
 $Y_{t_0}$ is Seshadri  $A_{t_0}$-ample {\rm (}respectively Seshadri $A_{t_0}$-big{\rm )}. Then there exists
 an open set $U$ in $B$, $t_0\in U$, such that $Y_t$ is Seshadri $A_t$-ample {\rm (}respectively Seshadri
 $A_t$-big{\rm )} for any $t\in U.$\end{theorem}

\begin{exs*}\label{cd1} Let us explicitly compute the cohomolgical dimension in some of the examples discussed in Section \ref{exampl}.

i) For a surface $Y$ in a polarized manifold $(X,A)$ which is Seshadri $A$-ample,  Theorem \ref{Big2}$(2)$,  asserts that $k-3\leq\cd(X\setminus Y)\leq k-2$.
For example, let  $Y$ be an irregular surface of $X=\mathbb P^k$ which is $A$-ample. E.g., take  $Y=C\times\mathbb P^1$ as in Proposition \ref{elliptic} and apply  Corollary \ref{Same0} above. Then a result of Hartshorne (see \cite[Theorem  8.5]{CDAV})  implies that $\cd(\mathbb P^k\setminus Y)=k-2$. On the other hand, by combining Proposition \ref{F1} and Corollary \ref{Same0}, we know that  the rational normal scroll $\mathbb F_e$ in $\mathbb P^{e+3}$, with $e=1,2,3,4$, is $\mathscr O_{\mathbb P^{e+3}}(1)$-ample. By a result of Verdi \cite{Verdi} (see also \cite[Corollary 4.2]{BV}), $\mathbb F_e$ is a set-theoretic complete intersection in $\mathbb P^{e+3}$. This implies that $\mathbb P^{e+3}\setminus\mathbb F_e$ is covered by $k-2=e+1$ affine open subsets of $\mathbb P^{e+3}$, whence $\cd(\mathbb P^{e+3}\setminus\mathbb F_e)\leq e=k-3$. Then  $\cd(\mathbb P^{e+3}\setminus\mathbb F_e)=e$ by statement $2)$ of Theorem \ref{Big2}. This shows that the bounds given in Theorem \ref{Big2}$(2)$,  are in general the best possible.

\smallskip

ii) Let us compute the cohomological dimension of the complement of $\mathbb P^2\times\mathbb P^1$ in $\mathbb P^5$ via the Segre embedding. By Proposition \ref{product}  and Corollary \ref{Same0}$(2)$, $\mathbb P^2\times\mathbb P^1$ is $\mathscr O_{\mathbb P^5}(1)$-ample. Therefore Theorem \ref{Big2}$(2)$,  predicts that 
$1\leq\cd(\mathbb P^5\setminus(\mathbb P^2\times\mathbb P^1))\leq 3$.  If $\cd(\mathbb P^5\setminus(\mathbb P^2\times\mathbb P^1))=1$ then  by \cite[Corollary 2.4, p. 229]{Ample} it would follow that 
$$H^2((\mathbb P^5)^{\an},\mathbb C)\cong H^2((\mathbb P^2\times\mathbb P^1)^{\an},\mathbb C),$$ 
and this is impossible because the first space is $1$-dimensional while the second  is $2$-dimensional. On the other hand, as it  is well known, $\mathbb P^2\times\mathbb P^1$ is given in $\mathbb P^5$ by the vanishing of the $(2\times 2)$-minors of the matrix
$$\begin{pmatrix}
x_0 & x_2 & x_4    \\
x_1 & x_3 & x_5\\
\end{pmatrix}.$$
Hence $\mathbb P^5\setminus(\mathbb P^2\times\mathbb P^1)$ is covered by three affine open subsets, so  that $\cd(\mathbb P^5\setminus(\mathbb P^2\times\mathbb P^1))\leq 2$. We thus conclude  that $\cd(\mathbb P^5\setminus(\mathbb P^2\times\mathbb P^1))= 2$.

\smallskip

iii) In the case of the Segre embedding $i\colon\mathbb P^3\times\mathbb P^1\hookrightarrow\mathbb P^7$, from Proposition \ref{product}  and Corollary \ref{Same0}$(2)$,  it follows that  $\mathbb P^3\times\mathbb P^1$ is $\mathscr O_{\mathbb P^7}(1)$-ample. Therefore Theorem \ref{Big2}$(2)$,  predicts that 
$2\leq\cd(\mathbb P^7\setminus(\mathbb P^3\times\mathbb P^1))\leq 5$. If $\cd(\mathbb P^7\setminus(\mathbb P^3\times\mathbb P^1))\leq 3$, using again \cite[Corollary 2.4, p. 229]{Ample} we deduce that $$H^2((\mathbb P^7)^{\an},\mathbb C)\cong H^2((\mathbb P^3\times\mathbb P^1)^{\an},\mathbb C),$$ and this is again absurd. Therefore $4\leq\cd(\mathbb P^7\setminus(\mathbb P^3\times\mathbb P^1))\leq 5$. On the other hand, as $\mathbb P^3\times\mathbb P^1$ is a rational normal scroll in $\mathbb P^7$, by 
\cite[Theorem 4.1]{BV}, $\mathbb P^3\times\mathbb P^1$ is defined  set-theoretically by five equations in $\mathbb P^7$. This implies that $\cd(\mathbb P^7\setminus(\mathbb P^3\times\mathbb P^1))\leq 4$, whence we finally get $\cd(\mathbb P^7\setminus(\mathbb P^3\times\mathbb P^1))=4$.
\end{exs*}

\begin{rem*} In Section \ref{codimone} we have seen that the normal bundle of a Seshadri $A$-big submanifold of codimension one in a polarized manifold of dimension $\geq 3$ need not be ample. However,  Paoletti proved in \cite{Pao} that the normal bundle of every Seshadri $A$-big smooth curve in a polarized threefold $(X,A)$ is ample. It would be interesting to see whether this result can be somehow generalized. For instance, is the normal bundle of every Seshadri $A$-big smooth curve in a polarized manifold $(X,A)$ of arbitrary dimension $\geq 4$ always ample?
\end{rem*}

\bigskip
\bigskip

 \begin{tabular}{l} 
 Lucian B\u adescu and Mauro C. Beltrametti\\
  Universit\`a degli Studi di Genova\\
    Dipartimento di Matematica\\
 Via Dodecaneso 35\\
  I-16146 Genova, Italy\\
E-mails: badescu@dima.unige.it and  beltrame@dima.unige.it \\
 \end{tabular}


\small

\begin{thebibliography}{999}

 \bibitem{AK} A. Altman and S. Kleiman,  \textit{Introduction to Grothendieck
 Duality Theory},
 Lecture Notes in Math. {\bf 146}, (1970), Springer-Verlag, New York, 1970.

\bibitem{Ba} L. B\u{a}descu,  Infinitesimal deformations of negative
weights and hyperplane sections, in {\em  Algebraic Geometry}, Proceedings
L'Aquila 1988, Lecture Notes in Math. {\bf 1417} (1990), Springer-Verlag, 1--22.

\bibitem{B1} L. B\u adescu, \textit{Algebraic Surfaces}, Universitext (translated by V. Mas\c ek), Springer-Verlag,
New York-Berlin-Heidelberg, 2001.

\bibitem{B2} L. B\u adescu, \textit{Projective Geometry and Formal 
Geometry}, Monografie Matematyczne, Instytut Matematyczny PAN, vol. 65, New Series,  Birkh\"auser, 2004.


\bibitem{BBF} L. B\u{a}descu, M.C. Beltrametti and P. Francia, Positive
curves in polarized manifolds, Manuscripta Math. {\bf  92} (1997), 369--388.

\bibitem{BV} L. B\u adescu and G. Valla, Grothendieck--Lefschetz theory, set-theoretic complete intersections and rational normal scrolls, Journal of Algebra  {\bf 324} (2010), 1636--1655.


\bibitem{BHPV} W. Barth, K. Hulek, C.A.M. Peters, and A. Van de Ven, \textit{Compact Complex Surfaces}, Ergebnisse der Math. vol. 4, Springer-Verlag, Berlin, 2004.

\bibitem{BS1} M.C. Beltrametti and A.J. Sommese, Remarks on numerically positive  and
big line bundles, in {\em  Projective Geometry
with Applications}, (ed. by E. Ballico), Lecture Notes in Pure and Applied Math. {\bf 166} (1994),
9--18, M. Dekker, New York.

\bibitem{BS} M.C. Beltrametti and A.J. Sommese, Notes on embeddings of
blowups, Journal of Algebra {\bf 186} (1996), 861--871. 


 \bibitem{BG} S. Bloch and D. Gieseker, The positivity of the Chern classes of
 an ample vector bundle, Invent. Math. {\bf 12} (1971), 112--117. 

 \bibitem{Co} M. Coppens, Embedding of blowing-ups, Universit\`a degli
 Studi di Bologna,
 Dipartimento di Matematica, Seminari di Geometria 1991-1993, 89--100.

\bibitem{Fu} W. Fulton, {\em Intersection Theory}, Ergebnisse der Math. Grenzgeb. (3) 2, Springer-Verlag, Berlin 1984.

\bibitem{FH} W. Fulton and J. Hanssen, A connectedness theorem for proper varieties with applications to intersections and singularities, Annals of Math. \textbf{110} (1979), 
159--166.

\bibitem{FL1} W. Fulton and R. Lazarsfeld, On the connectedness of degeneracy loci and special divisors, Acta Math. \textbf{146} (1981), 271--283.

\bibitem{FL} W. Fulton and R. Lazarsfeld, Positivity and excess
intersection, in {\em Enumerative and Classical Algebraic Geometry},
Nice 1981,  Progr.  Math. {\bf 24},  Birkhauser (1982), 97--105.

\bibitem{GH} Ph. Griffiths and J. Harris, \textit{Principles of Algebraic Geometry}, John Wiley \& Sons, New York, Chichester, Brisbane, Toronto, 1978.

 \bibitem{EGA} A. Grothendieck, {\em \'El\'ements de G\'eom\'etrie
Alg\'ebrique}, II, III, Premi\`ere Partie, Publ. Math. IHES 8, 11 (1961), 20 (1964).

\bibitem{Local} A. Grothendieck, {\em Local Cohomology}, Lecture Notes in
Math. {\bf  41}, Springer-Verlag, NewYork, 1967.

\bibitem{SGA} A. Grothendieck, {\em Rev\^{e}tements \'Etales et Groupe
Fondamental}, (SGA I), Lectures Notes in Math. {\bf 224}, Springer-Verlag, New York, 1971.

\bibitem{AVB} R. Hartshorne, Ample vector bundles,  Publ. Math. Inst. Hautes \'Etudes Sci.  no. 29 (1966), 63--94.

\bibitem{CDAV} R. Hartshorne, Cohomological dimension of algebraic
varieties, Ann. of Math. {\bf 88} (1968), 403--450.

\bibitem{Ample} R. Hartshorne, {\em Ample Subvarieties of Ample Varieties},
Lecture Notes in Math. {\bf 156}, Springer-Verlag, New York, 1970.

\bibitem{Hartshorne} R. Hartshorne, {\em Algebraic Geometry}, GTM {\bf 52}, Springer-Verlag, New York, 1977.

\bibitem{HM} H. Hironaka and H. Matsumura, Formal functions and formal
embeddings, J. Math. Soc. Japan {\bf 20} (1986), 52--82.

\bibitem{Kl} S.L. Kleiman, Ample vector bundles on algebraic surfaces,
Proc. Amer. Math. Soc. {\bf 21} (1969), 673--676.

\bibitem{Laz} R. Lazarsfeld, Some applications of the theory of positive vector bundles,
in {\em Complete Intersections} (Acireale, 1983, ed. by S. Greco and R. Strano), Lecture Notes in Math.  {\bf  1092}, 
Springer-Verlag (1984), 21--61.

\bibitem{Laz1} R. Lazarsfeld,  \textit{Positivity in Algebraic Geometry}, vols. I, II, Springer-Verlag, Berlin-Heidelberg-New York,  2004.

\bibitem{Pao}  R. Paoletti,  Seshadri positive curves in a smooth projective
 $3$-fold,  Atti Accad. Naz. Lincei Rend., Matematica e Applicazioni, Ser.
IX, vol. VI,  (1995), 259--274.

\bibitem{Pao2}  R. Paoletti,  Seshadri constants, gonality of space  curves and restriction of stable bundles, J. Differential Geometry {\bf 40} (1994), 475--504.

\bibitem{SW} A.K. Singh and U. Walther,  On the arithmetic rank of certain Segre products, in 
{\em Commutative Algebra and Algebraic Geometry},  Contemp. Math. {\bf 390}, Amer. Math. Soc., Providence, RI, 2005, pp. 147--155.

\bibitem{So} A.J. Sommese, Submanifolds of abelian varieties,  Math.
 Ann. {\bf 233} (1978), 229--256.


\bibitem{Sp} R. Speiser, Cohomological dimension and abelian varieties, Amer. J. Math. \textbf{95} (1973), 1--34.

\bibitem{Verdi} L. Verdi, Esempi di superficie e curve intersezioni complete insiemistiche, Boll. Unione Mat. Ital. (6) {\bf 5-A} (1986), 47--53.

\end{thebibliography}
\end{document}